\documentclass{article}
\usepackage[T1]{fontenc} % enables bold small caps (needed for \textproc/\textsc in section titles)

\usepackage[margin=.9in]{geometry}
\usepackage{graphicx}
\usepackage{amsmath,amssymb,bm}
\usepackage[dvipsnames]{xcolor}
\usepackage{float}
\usepackage{framed}
\usepackage{enumitem}
\usepackage{setspace}
\usepackage{booktabs}
\usepackage{soul}
\usepackage{cite}
\usepackage{pdflscape}
\usepackage{longtable}
\usepackage{adjustbox}

\usepackage{algorithm}
\usepackage[indLines=true,noEnd=false]{algpseudocodex}

\usepackage[
    colorlinks=true,
    linkcolor=Maroon,
    urlcolor=blue,
    citecolor=blue,
    anchorcolor=blue
]{hyperref}

\algrenewcommand\algorithmicrequire{\textbf{Input:}}
\algrenewcommand\algorithmicensure{\textbf{Output:}}

\definecolor{algstep}{RGB}{40,40,40}
\definecolor{algemph}{RGB}{0,76,153}

\newcommand{\R}{\mathbb{R}}

\newcommand{\cali}{\mathrm{i}}
\newcommand{\calI}{\mathcal{I}}
\newcommand{\lam}{\lambda}
\newcommand{\Lc}{{\cal L}_c}

\newcommand{\vecu}{\mathbf{u}}
\newcommand{\psis}{\psi^{(\mathrm{s})}}
\newcommand{\psin}{\psi^{(\mathrm{n})}}
\newcommand{\splus}{\scalebox{0.55}{$+$}}
\newcommand{\qedflush}{\hfill \scalebox{0.7}{$\blacksquare$}}

\newcommand{\PABIPP}{\textcolor{algemph}{\textproc{P-ABIPP}}}

\newcommand{\PAADMM}{\textcolor{algemph}{\textproc{P-Adapt-ADMM}}}
\newcommand{\AADMM}{\textcolor{algemph}{\textproc{Adapt-ADMM}}}
\newcommand{\JAADMM}{\textcolor{algemph}{\textproc{J-Adapt-ADMM}}}

\usepackage{amsmath}
\usepackage{amsfonts}
\usepackage{latexsym}
\usepackage{amssymb}
\usepackage{framed}

\newtheorem{thm}{Theorem}[section]
\newtheorem{theorem}[thm]{Theorem}

\newtheorem{lemma}[thm]{Lemma}
\newtheorem{prop}[thm]{Proposition}
\newtheorem{proposition}[thm]{Proposition}
\newtheorem{cor}[thm]{Corollary}

\newtheorem{definition}[thm]{Definition}

\newtheorem{remark}[thm]{Remark}

\newcommand{\beq}{\begin{equation}}
\newcommand{\eeq}{\end{equation}}
\newcommand{\beqa}{\begin{eqnarray}}
\newcommand{\eeqa}{\end{eqnarray}}
\newcommand{\beqas}{\begin{eqnarray*}}
\newcommand{\eeqas}{\end{eqnarray*}}
\newcommand{\bi}{\begin{itemize}}
\newcommand{\ei}{\end{itemize}}

\newcommand{\vgap}{\vspace{.1in}}

\newcommand{\Ik}{{\mathcal{I}}_k}

\newcommand{\inner}[2]{\langle #1,#2\rangle}
\newcommand{\Inner}[2]{\left \langle #1\,,#2 \right \rangle}

\newcommand{\interior}{\mathrm{int}\,}

\newcommand{\dom}{\mathrm{dom}\,}

\newcommand{\bConv}[1]{\overline{\mbox{\rm Conv}}\,(\R^{#1})}

\newcommand{\ty}{\tilde y}
\newcommand{\tT}{\tilde T}

\newcommand{\tq}{\tilde q}

\title{A Parallel Adaptive Proximal ADMM for Linearly Constrained Separable Composite Optimization Problems}
\begin{document}
\author{Leandro F. Maia\thanks{School of Mechanical, Industrial, and Manufacturing Engineering, Oregon State University,
Corvallis, OR, 97331. \protect\protect\href{mailto:leandro.maia@oregonstate.edu}{fariasmaia@gmail.com}}\hspace*{0.5em} \and Renato D.C. Monteiro\thanks{School of Industrial and Systems Engineering, Georgia Institute of
Technology, Atlanta, GA, 30332-0205. 
\protect\protect\href{mailto:monteiro@isye.gatech.edu}{monteiro@isye.gatech.edu}}\and Gilson N. 
Silva\thanks{Department of Mathematics, Federal University of Piauí, Teresina, PI, 64049-550. \protect\protect\href{mailto:monteiro@isye.gatech.edu}{gilson.silva@ufpi.edu.br}}}

\maketitle

\begin{abstract}
We propose and analyze an adaptive proximal alternating direction method of multipliers (ADMM), based on Jacobi-type updates, for linearly constrained separable nonconvex composite optimization problems. The objective consists of two block-separable functions, one continuously differentiable and weakly convex, and the other convex and Lipschitz continuous. The proposed method is adaptive to all problem parameters and allows inexact solutions of the block proximal subproblems. At each iteration, we solve all block subproblems in parallel, then apply a novel test to decide whether to update the Lagrange multiplier and/or the penalty parameter. Without any rank assumptions on the constraint matrices, we establish an ${\cal O}(\epsilon^{-3})$ iteration-complexity bound for obtaining an approximate first-order stationary point. This matches the best-known complexity for cyclic proximal ADMM methods and improves the existing complexity of Jacobi-type schemes in this general nonconvex setting. Numerical experiments on a nonconvex quadratic program illustrate the proposed method's computational benefits.

\vgap


\end{abstract}

\section{Introduction}
In this article, we provide the iteration-complexity analysis for the Parallel Adaptive Proximal Alternating Direction Method of Multipliers, referred to as {\PAADMM}, when solving the linearly constrained, composite optimization problem  
\begin{equation}\label{initial.problem}
    \phi^* = \min_{x\in\R^n} \left\{ \phi(x) := f(x)+h(x)  :   Ax=b \right\},  
\end{equation}
where $h:\R^n \to (-\infty,\infty]$ is a proper closed convex function which is Lipschitz continuous on its compact domain, $f:\R^{n} \to \R$ is a continuously differentiable, weakly convex function on the domain of $h$, $A:\R^{n} \to \R^l$ is a linear operator, and $b\in \R^l$ is a vector in the image of $A$. Moreover, for some positive integer $B$ (the number of blocks), with $B\geq 2$, and positive integer vector $(n_1,\ldots,n_B)$ such that $n=\sum_{t=1}^{B} n_{t}$, the functions $f(\cdot)$ and $h(\cdot)$, and the linear operator $A(\cdot)$ have the following blockwise representation
\begin{gather}
\begin{gathered}
f(x)=\sum_{t=1}^Bf_t(x_{t}), \quad h(x)=\sum_{t=1}^{B}h_{t}(x_{t}), \quad Ax=\sum_{t=1}^{B}A_{t}x_{t},
\end{gathered}
\label{eq:block_structure}
\end{gather}
where, for each $t\in\{1,\ldots,B\}$, $h_t:\mathbb{R}^{n_t}\to(-\infty,\infty]$
is a proper closed convex function with compact domain, $f_{t}:\R^{n_t}\to \R$ is a continuously differentiable, weakly convex function on the domain of $h_t$, and $A_{t}:\R^{n_t} \to \R^l$ is a linear operator. 

ADMM is an efficient block-decomposable first-order method for solving \eqref{initial.problem}--\eqref{eq:block_structure} (see e.g., the survey \cite{BoydEtAl2011}). Two common strategies for updating the primal blocks are the Gauss--Seidel and Jacobi schemes, often referred to as cyclic and parallel updates, respectively. In the Gauss--Seidel scheme, the blocks are updated sequentially, each update using the latest values of the preceding blocks. In the Jacobi scheme, all blocks are updated from the information available at the start of the iteration, so the block subproblems can be solved simultaneously, and the new primal point is formed after all of them are completed. Their applications include image restoration, semidefinite programming, robust principal component analysis, phase retrieval, and matrix completion for the former (see, e.g., \cite{GoldsteinOsher2009,WenGoldfarbYin2010,TaoYuan2011,WenYangLiuMarchesini2012,XuYinWenZhang2012}), and exchange and consensus problems, graphical model selection, and multi-period optimal power flow for the latter (see, e.g., \cite{DengLaiPengYin2017,ChangHongWang2015,HeMaXuYuan2022,SubramanyamEtAl2021}).
%A formal description of the parallel proximal ADMM framework considered in this paper is given in Subsection~\ref{sub:parallel.prox.ADMMM}.

For a given tolerance pair $(\rho,\eta)\in \mathbb{R}_{++}^2$, our goal is to compute a quadruple $(x,p,v,\varepsilon)\in \dom h \times A(\R^n) \times \R^n \times \R_+$ satisfying
\begin{equation}\label{eq:rho.eta.solution}
 v \in \nabla f(x) + \partial_\varepsilon h(x) + A^*p, \quad  \|v\|^2 + \varepsilon \le \rho^2, \quad \|Ax-b\|\leq \eta.
\end{equation}
We call such a quadruple a $(\rho,\eta)$-stationary solution of \eqref{initial.problem}--\eqref{eq:block_structure} throughout the paper. For later reference, we call a quadruple satisfying the first two conditions in \eqref{eq:rho.eta.solution} a $\rho$-stationary solution, and a point $x$ satisfying the last condition $\eta$-feasible.

\paragraph{Contributions.}
We show  that {\PAADMM} obtains a $(\rho,\eta)$-stationary solution of \eqref{initial.problem}--\eqref{eq:block_structure} in at most ${\cal O}(B^2\rho^{-2}\eta^{-1})$ iterations, and hence in ${\cal O}(\epsilon^{-3})$ iterations where $\epsilon:=\min\{\rho,\eta\}$, under suitable choices of its input parameters. This bound is established without assuming any rank condition on the constraint matrices or the last-block condition \eqref{cond:last-block} (see the ``Related work'' discussion later in this section). {\PAADMM} also has the following properties:
\begin{itemize}
\item {\it Adaptiveness:} It is fully adaptive, i.e., it requires no knowledge of the weak convexity parameters $m_t$'s of the $f_t$'s or of the Lipschitz constant of $h$. In particular, the block prox stepsizes are chosen by a backtracking procedure whose acceptance test involves only quantities available to the corresponding block, and hence can be carried out by all blocks in parallel.

\item {\it Inexactness:} It allows each block proximal subproblem to be solved either exactly or inexactly, according to a relative stationarity criterion that can be verified by an adaptive accelerated first-order method.

\item {\it Novel multiplier update test:} The multiplier update test replaces the fixed thresholds used in \cite{MaiaEtAl2026}, allowing more frequent multiplier updates.
\end{itemize}

\paragraph{Related works.}
Since ADMM and its proximal variants have been extensively studied in both the convex (see, e.g., \cite{Gabay1983,EcksteinBertsekas1992,MonteiroSvaiter2013,BoydEtAl2011}) and the nonconvex (see, e.g., \cite{HongLuoRazaviyayn2016,WangYinZeng2019,KongMeloMonteiro2023,KongMonteiro2024,MaiaEtAl2026,phdMaia2025}) settings, we focus on methods based on the Jacobi schemes. We first discuss convex and nonconvex ADMMs with Jacobi block updates, and we then compare {\PAADMM} with its Gauss--Seidel counterpart \cite{MaiaEtAl2026}.

Jacobi-type schemes for convex problems have a long history (see e.g.,
\cite{ChenTeboulle1994,Ruszczynski1995,HanYuanZhang2014,
ChatzipanagiotisDentchevaZavlanos2015,LinLiuLi2015,WangHongMaLuo2015,
HeYuan2015blockwise,He2009parallel,HeTaoYuan2015,Tao2014parallel,HeHouYuan2015,DengLaiPengYin2017,WangBanerjeeLuo2014}). In particular, it is shown by example in \cite{HeHouYuan2015} that the
direct Jacobi extension of the augmented Lagrangian method may diverge even
in the convex case, and restored global convergence, with an
$\mathcal O(1/k)$ rate, by adding a relaxation step whose stepsize depends
on the number of blocks, under a full-column-rank assumption on the
$A_t$'s. In \cite{DengLaiPengYin2017} instead, it is added proximal
terms to the block subproblems, and established global convergence and an $o(1/k)$ rate without any rank assumption. Related parallel schemes for convex problems
include the linearized ADMM with parallel splitting and adaptive penalty of
\cite{LinLiuLi2015}, the two-block reformulation of
\cite{WangHongMaLuo2015}, and the randomized parallel direction method of
multipliers \cite{WangBanerjeeLuo2014}. All of these works assume the objective is convex, and their rates are given either for ergodic
objective-gap and feasibility measures, or for the distance between consecutive iterates, instead of the stationarity residual pair in \eqref{eq:rho.eta.solution}. In addition, their block subproblems are solved
exactly or with a single linearized step, and their relaxation, proximal,
or penalty parameters are chosen based on the problem data.

Results for nonconvex problems are more recent (see e.g., \cite{LiuShenGu2019,MeloMonteiro2017,ZhangMaZhang2020,SubramanyamEtAl2021,YangEtAl2022,DuEtAl2026}) and, to a large extent, rely on a special last block. More specifically, one block, say block $B$, is unconstrained (i.e., $h_B\equiv 0$), has a differentiable objective with Lipschitz continuous gradient, and its constraint matrix satisfies
\begin{equation}
 \operatorname{Im}(A_B)
 \supseteq
 \{b\}\cup\operatorname{Im}(A_1)\cup\cdots\cup
 \operatorname{Im}(A_{B-1}).
 \tag{LB}\label{cond:last-block}
\end{equation}
We refer to these requirements as the \emph{last-block condition}. It is worth noting that the choice $A_B=I$ automatically satisfies the range part of this condition. It is established in \cite{MeloMonteiro2017} that an $\mathcal O(\epsilon^{-2})$ iteration-complexity bound holds for a nonconvex Jacobi-type proximal ADMM under \eqref{cond:last-block}, with block subproblems solved exactly and prox and penalty parameters chosen from problem constants. In \cite{LiuShenGu2019}, a parallel linearized ADMM is studied under a version of \eqref{cond:last-block} together with a full-column-rank assumption, and only asymptotic convergence is obtained.

Two works dispense \eqref{cond:last-block} by applying an ADMM-type method to a penalty reformulation of \eqref{initial.problem}--\eqref{eq:block_structure} that artificially satisfies it. In \cite{ZhangMaZhang2020} it is obtained an $O(\epsilon^{-4})$ iteration complexity. In \cite{SubramanyamEtAl2021} it is proposed an adaptive parallel Jacobi scheme based on the same mechanism and also obtained an $\mathcal O(\epsilon^{-4})$ bound, under the additional assumptions that $[A_1,\ldots,A_B]$ has full row rank and that every block update returns a local minimizer. Finally, two recent works \cite{YangEtAl2022,DuEtAl2026} avoid \eqref{cond:last-block} without reformulating the problem, but they do not provide iteration-complexity bounds. In summary, for Jacobi-type ADMMs that do not assume the last-block condition, the best iteration-complexity bound we are aware of is $\mathcal O(\epsilon^{-4})$ \cite{ZhangMaZhang2020,SubramanyamEtAl2021}.

We end this discussion by comparing the current work with that in \cite{MaiaEtAl2026}, which applies to a nonseparable version of \eqref{initial.problem}--\eqref{eq:block_structure}. It also allows inexact block solves, adapts to unknown weak convexity and Lipschitz constants, and requires neither a rank condition on $A$ nor the last-block condition \eqref{cond:last-block}. However, its block updates follow a sequential Gauss--Seidel order. By adopting a parallel block-update scheme, the convergence analysis introduces additional terms that are not present in the Gauss--Seidel scheme. In terms of tolerance only, both iteration complexities for finding a stationary solution are $\tilde{\cal O}(\epsilon^{-3})$, where $\epsilon=\min\{\eta,\rho\}$. To the best of our knowledge, {\PAADMM} is the first Jacobi-type proximal ADMM with an $\tilde{\cal O}(\epsilon^{-3})$ iteration-complexity bound that requires neither the last-block condition nor a rank condition on $A$.

\paragraph{Structure of the Paper.} The remainder of the paper is organized as follows. Subsection~\ref{sub:notation} presents the notation, definitions, and basic facts used throughout the paper. Section~\ref{sec:prelim} provides an overview of {\PAADMM} and details the notion of inexact solution of the block proximal subproblems adopted in this work. Section~\ref{sec:algorithm} states the assumptions imposed on problem~\eqref{initial.problem}--\eqref{eq:block_structure}, presents {\PABIPP}, the subroutine that performs a parallel adaptive block inexact proximal point iteration, and states its main properties. It also details {\PAADMM} and states its iteration-complexity bound. Sections~\ref{sec:proof.prop} and~\ref{sec:proof.thm} are technical sections dedicated to proving the main results of {\PABIPP} and {\PAADMM}, respectively. Section~\ref{sec:numerics} presents proof-of-concept numerical experiments on a class of nonconvex linearly constrained quadratic problems. Section~\ref{sec:conclusion} presents some concluding remarks. Appendix~\ref{app:technical} collects technical results used in the paper, while Appendix~\ref{app:adaptive.tau} describes a variant of {\PABIPP} with an adaptively chosen relaxation parameter.

\subsection{Notation, Definitions, and Basic Facts}\label{sub:notation}

This subsection introduces the fundamental notation employed throughout the paper. The set of real numbers is denoted by $\R$, while $\R_{+}$ and $\R_{++}$ represent the sets of non-negative and positive real numbers, respectively.
The $n$-dimensional Euclidean space $\R^n$ is equipped with an inner product $\Inner{\cdot}{\cdot}$.

The norm induced by $\Inner{\cdot}{\cdot}$ is denoted by $\|\cdot\|$.
$\R^n_{++}$ and $\R^n_{\splus}$ denote the sets of vectors in $\R^n$ with positive and non-negative entries, respectively.
The smallest positive singular value of a nonzero linear operator $Q:\R^n\to \R^l$ is denoted by $\nu^{\splus}_Q$. The operator norm is defined as $\|Q\|:=\sup \{\|Q(w)\|:\|w\|=1\}$.
For a symmetric and positive definite matrix $S$, the norm induced by $S$ on $\R^n$ is defined as $\|\cdot\|_S:=\Inner{\cdot}{S(\cdot)}^{1/2}$.
Given $x=(x_{1},\ldots,x_{B}) \in \mathbb{R}^{n_{1}} \times \cdots \times \mathbb{R}^{n_{B}}$, the following aggregated quantities are defined:
\begin{gather*}
\begin{gathered}
x_{<t}:=(x_{1},\ldots,x_{t-1}), \quad  x_{>t}:=(x_{t+1},\ldots,x_{ B}),\quad x_{\leq t}:=(x_{<t},x_{t}),\quad x_{\geq t}:=(x_{t},x_{>t}).
\end{gathered}
\end{gather*}
Moreover, given a differentiable function $g:\R^{n_1}\times\cdots\times\R^{n_B}\to\R$ and $t\in\{1,\ldots,B\}$, $\nabla_t g(x)$ denotes the partial gradient of $g$ with respect to the block $x_t$ at $x$.
For $a\in \R^B$, define
\begin{equation}\label{def.lambda.min.max}
\min(a)  = \min_{1\leq t\leq B} a_t  \quad \text{and} \quad \max(a)  = \max_{1\leq t\leq B} a_t .
\end{equation}
For a closed convex set $Z \subset \R^n$, $\partial Z$ denotes its boundary.
The distance from a point $z \in \R^n$ to $Z$, measured with $\|\cdot\|$, is denoted by ${\rm dist}(z,Z)$.
The indicator function of $Z$, denoted $\delta_Z$, is defined by $\delta_Z(z)=0$ if $z \in Z$ and $\delta_Z(z)=+\infty$ otherwise.
Let $\log(\cdot)$ denote the base-2 logarithm.

For a function $g :\R^n\to (-\infty,\infty]$, the effective domain is denoted by $\dom g := \{x\in \R^n : g(x) < +\infty\}$.
A function $g$ is proper if $\dom g \ne \emptyset$. The set of all lower semi-continuous proper convex
functions defined on $\R^n$ is denoted by $\bConv{n}$.
For $\varepsilon\geq 0$, the $\varepsilon$-subdifferential of $g\in \bConv{n}$ at $x\in \dom g$ is defined as
\begin{equation*}
\partial_\varepsilon g(x):=\{ u\in \R^n ~:~ g(\tilde x)\geq g(x)+\Inner{u}{\tilde x-x}-\varepsilon, \forall \tilde x\in \R^n\}.
\end{equation*}
When $\varepsilon=0$, the $\varepsilon$-subdifferential coincides with the classical subdifferential: $\partial g(\cdot):=\partial_0 g(\cdot)$.
It is well-known (see \cite[Prop.\ 1.3.1 of Ch.\ XI]{lemarechal1993})
that for any $\beta>0$ and  $g \in \bConv{n}$,
\begin{equation}\label{prop1:subd}
\partial_{\varepsilon} (\beta g)(\cdot ) = \beta \partial_{(\varepsilon/\beta)} g(\cdot ).
\end{equation}
Furthermore, if $ h_t \in \bConv{n_t}$ for $t=1,\ldots,B$ and $h(x):=\sum_{t=1}^Bh_t(x_t)$ for any $x=(x_1,\ldots,x_B)$, then
the following property holds:
(see \cite[Remark 3.1.5 of Ch.\ XI]{lemarechal1993})
\begin{align}
\label{prop2:subd}
\partial_{\varepsilon}h(x) =\bigcup \Big\{ \partial_{\varepsilon_1}h_1(x_1) \times \cdots \times
\partial_{\varepsilon_B}h_B(x_B) : \varepsilon_t\ge 0,\ \varepsilon_1+\cdots+\varepsilon_B\le \varepsilon\Big\}.
\end{align}

\section{Preliminaries}\label{sec:prelim}

This section is organized into two subsections. Subsection~\ref{sub:overview} presents an overview of {\PAADMM}, and Subsection~\ref{sub:inexact} introduces the notion of the inexact solution used for the block proximal subproblems.

\subsection{An overview of {\PAADMM}}
\label{sub:overview}

In this subsection, we outline the main algorithm of this paper and defer the details for later. {\PAADMM} is based on the augmented Lagrangian function
\begin{equation}\label{DP:AL:F}
{\cal L}_c(x;p)
:=
\phi(x)+\langle p,A x-b\rangle+\frac{c}{2}\|A x-b\|^2,
\end{equation}
where $p$ is the Lagrange multiplier and $c>0$ is the penalty
parameter.

We next outline an iteration of {\PAADMM}, which depends on a test for deciding when to update the Lagrange multiplier $p$. Since the details of this \textit{Lagrange multiplier test} (LMT) are non-trivial, we postpone them to Subsection~\ref{sub:PAADMM} and exclude them from the outline below for simplicity. Let $\tau:=1/[4(B-1)]$, where $B\geq 2$ is the number of blocks.  Given $(x,p,c,\gamma)\in \dom h\times A(\R^n)\times\R_{++}\times\R_{++}^B$, the method either terminates or generates a new quadruple $(x^{\splus},p^{\splus},c^{\splus},\gamma^{\splus})$ as follows:

\begin{enumerate}
\item[(i)] In parallel, for every $t=1,\ldots,B$, find a pair $(\hat x^{\splus}_t,\hat \gamma^{\splus}_t)$ such that $\hat x_t^{\splus}$ is a near-stationary solution of
\begin{equation*}
\min_{u_t\in\R^{n_t}}\left\{\hat\gamma_t^{\splus}
{\cal L}_c(x_{<t},u_t,x_{>t};p)
+\frac12\|u_t-x_t\|^2
\right\}
\end{equation*}
and yields a sufficient decrease on the function $u_t \mapsto \hat\gamma_t^{\splus} {\cal L}_c(x_{<t},(1-\tau)x_t + \tau u_t ,x_{>t};p)$, i.e., the quantity
\[\hat\gamma_t^{\splus}
\left[{\cal L}_c(x;p)-{\cal L}_c
\bigl(x_{<t},(1-\tau)x_t+\tau \hat x_t^{\splus},x_{>t};p\bigr)
\right]
\]
is positive and bounded away from zero in a relative sense. Then, set $\hat x^{\splus}=
(\hat x_1^{\splus},\ldots,\hat x_B^{\splus})$ and $\gamma^{\splus}=
(\hat\gamma_1^{\splus},\ldots,\hat\gamma_B^{\splus})$.
\item[(ii)] If $\hat x^{\splus}$ is a near-stationary solution of $\min_{u\in\R^n}{\cal L}_c(u;p)$, then set
\[
x^{\splus} = \hat x^{\splus},
\qquad
p^{\splus}=p+c(A \hat x^{\splus}-b), \qquad c^{\splus}=2c.
\]
Moreover, if $x^{\splus}$ is  a near-feasible point of \eqref{initial.problem}, then terminate {\PAADMM};
in this case, $(x^{\splus},p^{\splus})$ is a near-stationary solution of \eqref{initial.problem}--\eqref{eq:block_structure} in the sense of \eqref{eq:rho.eta.solution}. 
Otherwise, if $\hat x^{\splus}$ is not a near-stationary solution of $\min_{u\in\R^n}{\cal L}_c(u;p)$, then set
\[
x^{\splus}=(1-\tau)x+\tau \hat x^{\splus},
\qquad
c^{\splus}=c,
\]
and
\[
p^{\splus}
=
\begin{cases}
p+c(A \hat x^{\splus}-b),
& \text{if LMT succeeds},\\
p,
& \text{otherwise}.
\end{cases}
\]
% The multiplier update test is a new rule that determines
% whether to accept the candidate multiplier.
\end{enumerate}

Unless {\PAADMM} terminates in (ii), the next iteration starts from $(x^{\splus},p^{\splus},c^{\splus},\gamma^{\splus})$. Section~\ref{sec:algorithm} gives the precise implementation of the outline. In particular, Subsection~\ref{sub:PABIPP} describes and studies the main properties of {\PABIPP}, a subroutine that implements step (i) above. Subsection~\ref{sub:PAADMM} describes {\PAADMM}, the main algorithm of this paper that implements steps (i) and (ii), and which, through multiple calls of {\PABIPP} and the LMT, builds a near-stationary solution of \eqref{initial.problem}--\eqref{eq:block_structure} in the sense of \eqref{eq:rho.eta.solution}. This subsection also states the iteration complexity of {\PAADMM}.

% \begin{itemize}
%     \item 
%     In parallel, for every $t=1,\ldots,B$, find a pair $(\hat z^+_t,\hat \lam^+_t)$ such that $\hat z_t^+$ is a near-stationary solution of ???? and yields a sufficient decrease on the function $\hat\lam_t^+ {\cal L}_c(z_{<t},\cdot(1-\tau)z + \tau \hat z^+,z_{>t};q)$, i.e., the quantity
%     \[
% \hat\lam_t^+ \left[ {\cal L}L_c(z_{<t},z_t,z_{>t};q) - {\cal L}L_c(z_{<t},(1-\tau)z_t + \tau \hat z^+_t,z_{>t};q) \right] 
%     \]
%     is sufficiently positive; finally, set
%     ????
%     \item 
%     test if $\hat z^+$ is a near-stationary point of (11); if so, then set
%     \[
%    (z^+, \lam^+)=(\hat z^+, \hat \lam^+), \qquad c^+=2 c, \qquad q^+=q+ c(A\hat z^+-b) 
%     \]
%     \item 
%     otherwise, set 
%     \[
%    (z^+, \lam^+)=((1-\tau)z + \tau \hat z^+, \hat \lam^+), \qquad c^+=c, \qquad q^+= \left\{
%    \begin{array}{ll} q+ c(A\hat z^+-b), & \mbox{if LMT succeeds} \\
% q & \mbox{otherwise}
%    \end{array} \right.
%     \]
%     % also, check whether to perform a Lagrange multiplier; if so set $q^+$ as above; otherwise set $q^+=q$ above;
% \end{itemize}

% outline...mentioning sections...

% Finally, the whole method terminates when $\|A\hat z^+-b\|$ is small

\subsection{Inexact Solution Criteria for the Block Subproblems}\label{sub:inexact}

This subsection details the solution criteria adopted for the block proximal subproblems of the {\PAADMM}. We consider the following composition optimization problem
\begin{align}\label{ISO:problem}
 \min \{\psi(u):=\psis(u)+\psin(u) : u\in\R^d\},   
\end{align}
where $\psis:\R^d \to \R$ is a differentiable function and $\psin\in \bConv{d}$ is proper, closed, and convex.

The following definition explains the inexact solution concept for the minimization problem \eqref{ISO:problem}.
\begin{definition}\label{def:rel.stationary}
For a given $y^0 \in \dom \psin$ and parameter
$\sigma \in \R_{\splus}$,
a triple $(z,w,\epsilon) \in \R^d\times\R^d\times\R_{\splus}$ satisfying
\begin{align}\label{ISO:Cond:1&2}
w\in\nabla\psis(z)
+\partial_{\epsilon}\psin(z),
\qquad
2\|w\|^2+\epsilon
\leq\sigma\|y^0-z\|^2 
\end{align}
is called a $(\sigma;y^0)$-relative stationary solution of \eqref{ISO:problem} with composite term $\psin$.
\end{definition}

We next specialize Definition~\ref{def:rel.stationary} to a block proximal subproblem. For $c>0$ and $q\in\R^l$, let
\begin{equation}\label{def:smooth:ALM}
{\cal L}_{c}^{(\mathrm{s})}(y;q)
:=
f(y)+\langle q,Ay-b\rangle+\frac{c}{2}\|Ay-b\|^2
\end{equation}
denote the differentiable part of the augmented Lagrangian ${\cal L}_{c}(y;q)$ defined in~\eqref{DP:AL:F}. Given a reference point $y\in\dom h$, a block
$t\in\{1,\ldots,B\}$, and a trial prox stepsize
$\lambda_t^{\splus}>0$, set
\begin{align}
\psis(u_t)=\lambda_t^{\splus}
{\cal L}_c^{(\mathrm{s})}(y_{<t},u_t,y_{>t};q)
+\frac12\|u_t-y_t\|^2, \qquad \psin(u_t)
=\lambda_t^{\splus}h_t(u_t),
\qquad u_t\in\R^{n_t}.
\end{align}
Definition~\ref{def:rel.stationary} is then applied with starting point $y^0=y_t$, a constant $\sigma>0$, and $(\psis,\psin)$ as above, yielding a $(\sigma;y_t)$-relative stationary solution $(z_t^{\splus},w_t^{\splus},\epsilon_t^{\splus})$ satisfying~\eqref{ISO:Cond:1&2}. Since all blocks use the same reference point $y$, all $B$ subproblems can be solved in parallel.

In general, computing an exact solution or a relative stationary point of \eqref{ISO:problem} is challenging. Proposition 3.5 of \cite{HeMonteiro2015} (see also Subsection 2.3 in \cite{fistaReport2021} and Appendix A in \cite{sujanani2023adaptive}) establishes the iteration-complexity of a variant of the original Nesterov's accelerated gradient method \cite{Nes83} (see also \cite{nesterov2013gradient} and \cite[Chapter 2]{nesterov1983method}) for obtaining an approximate solution, assuming that $\psis$ is either convex or strongly convex. Moreover, for $\sigma>0$, \cite[Appendix~B]{MaiaEtAl2026} describes an adaptive accelerated first-order method for computing such approximate solution. The variational interpretation of this inclusion in the smooth setting is also discussed in \cite[Appendix~C]{MaiaEtAl2026}, which it was related to directional derivatives and to nearby points with small exact subdifferential residuals.

\section{Algorithm and Main Results}\label{sec:algorithm}

This section presents the two steps of {\PAADMM} outlined in Subsection~\ref{sub:overview}. We first state the assumptions on problem~\eqref{initial.problem}--\eqref{eq:block_structure} that are used throughout the paper. Subsection~\ref{sub:PABIPP} then presents {\PABIPP}, the subroutine that implements step (i) of the outline, and states its main properties in Proposition~\ref{prop:PABIPP}. Subsection~\ref{sub:PAADMM} presents {\PAADMM}, which implements steps (i) and (ii), and states its iteration complexity in Theorem~\ref{thm:main}.

The main problem studied in this paper is problem~\eqref{initial.problem} with the block structure in~\eqref{eq:block_structure}. We assume that the vector $b\in\mathbb{R}^{l}$, the linear operator $A:\mathbb{R}^{n}\to\mathbb{R}^{l}$, and the block functions $f_t:\R^{n_t}\to\R$ and $h_t:\mathbb{R}^{n_t}\to(-\infty,\infty]$ for $t=1,\ldots,B$ satisfy the following conditions:

\begin{enumerate}

\item[(A1)] For every $t=1,\ldots,B$,
$h_t\in\bConv{n_t}$ is prox-friendly
(i.e., its proximal operator is easily computable) and has
compact domain ${\cal H}_t:=\dom h_t$.
Moreover, $h$ is $M_h$-Lipschitz continuous on
${\cal H}:={\cal H}_1\times\cdots\times{\cal H}_B$
for some $M_h\geq 0$;

\item[(A2)] For every $t=1,\ldots,B$, there exists $m_t\geq 0$
such that $f_t$ is $m_t$-weakly convex on ${\cal H}_t$, i.e.,
\[
f_t(\cdot)+\delta_{{\cal H}_t}(\cdot)
+\frac{m_t}{2}\|\cdot\|^2
\quad\text{is convex}.
\]
Moreover, $f_t$ is a continuously differentiable function on ${\cal H}_t$;

\item[(A3)] $A$ is nonzero and there exists $\bar{\mathrm{x}}\in \{x\in {\cal H}:Ax=b\}\neq \emptyset$ such that $\bar{d}:=\text{dist}(\bar{\mathrm{x}}, \partial {\cal H})>0$.  
\end{enumerate}

We now make some remarks about the above assumptions. First, since ${\cal H}$ is compact by (A1), the scalars
\begin{equation}\label{def:damH}
       D_h:=\sup_{z \in {\cal H}} \|z-\bar{\textrm{x}}\|, \quad \widehat \nabla_f:=\sup_{u \in {\cal H}} \| \nabla f(u)\|, \quad  \phi_{\mathrm{inf}} := \inf_{u \in {\cal H}} \phi(u), \quad  \phi_{\mathrm{sup}} := \sup_{u \in {\cal H}}  \phi(u),  
\end{equation}
are bounded. Second, assumption (A2) allows $m_t=0$, in which case $f_t$ is convex on ${\cal H}_t$.
Finally, it follows from (A3) that
$\mathrm{\bar x}$ is a Slater point for \eqref{initial.problem}
in the strong sense that $\mathrm{\bar x} \in \interior( {\cal H})$ and $A\mathrm{\bar x}=b$.

\subsection{The {\PABIPP}}\label{sub:PABIPP}
This subsection introduces \PABIPP, describes its relevant remarks, and states its main properties in Proposition~\ref{prop:PABIPP}. The formal description of {\PABIPP} is presented below.

{\floatname{algorithm}{Subroutine}
\renewcommand{\thealgorithm}{}
\begin{algorithm}[H]
\setstretch{1}
\caption{{\PABIPP}}
\begin{algorithmic}[1]
\Statex \hskip-1.8em
\Require $(y,q,\lam,c)
\in {\cal H}\times A(\R^n)\times\R^B_{++}\times\R_{++}$
\Ensure $(y^{\splus},\lam^{\splus};z^{\splus},r^{\splus},\delta^{\splus})
\in {\cal H}\times\R^B_{++}\times{\cal H}\times\R^n\times\R_{+}$ that satisfies the conclusions of Proposition~\ref{prop:PABIPP}.

\vspace{1em}

\Statex \hskip-2.8em \textbf{\textcolor{algemph}{STEP 1: Parallel adaptive block IPP iteration}}

\State $\bar\tau\gets 1/[4(B-1)]$
\label{ln:PABIPP:tau}

\For{$t=1,\ldots,B$ \textbf{in parallel}}
\label{ln:PABIPP:for}

    %\Statex
    \State $\lambda_t^{\splus}\gets\lambda_t$

    \State Compute a $\left(1/32;y_t\right)$-relative stationary solution
    $(z_t^{\splus},w_t^{\splus},\epsilon_t^{\splus})$ of
    \begin{equation}\label{algo.PABIPP-subproblem}
    \hspace{2.5cm}\min_{u\in\R^{n_t}}
    \left\{
    \lambda_t^{\splus}\Lc^{(\textrm{s})}(y_{<t},u_t,y_{>t};q)
    +
    \frac12\|u_t-y_t\|^2
    +
    \lambda_t^{\splus}h_t(u_t)
    \right\},\hspace{5.5em}
    \end{equation}
    with composite term $\lambda_t^{\splus}h_t(\cdot)$;
    see Definition~\ref{def:rel.stationary}.
    \label{ln:PABIPP:solve}

    \If{
\begin{equation}\label{algo:PABIPP.adaptive.lambda}
    \hspace{1cm}\Lc(y;q)
    -
    \Lc(y_{<t},y_t+\bar\tau(z_t^{\splus}-y_t),y_{>t};q)
    <
    \frac{\bar\tau^2}{2\lambda_t^{\splus}}\|z_t^{\splus}-y_t\|^2
    +
    \frac{c\bar\tau}{4}\|A_t(z_t^{\splus}-y_t)\|^2
\end{equation}
    }
    \label{ln:PABIPP:test}

        \State $\lambda_t^{\splus}\gets\lambda_t^{\splus}/2$
        and go to line~\ref{ln:PABIPP:solve}
        \label{ln:PABIPP:halve}

    \EndIf

    %\Statex

\EndFor
\label{ln:PABIPP:endfor}

%\Statex

\State $z^{\splus}\gets(z_1^{\splus},\ldots,z_B^{\splus})$,
$y^{\splus}\gets y+\bar\tau(z^{\splus}-y)$ and $\lambda^{\splus}\gets (\lambda^{\splus}_1,\ldots,\lambda^{\splus}_B)$
\label{ln:PABIPP:update}

\vspace{1em}

\Statex \textbf{\textcolor{algemph}{STEP 2: Computation of the residual pair $(r^{\splus},\delta^{\splus})$}}

\For{$t=1,\ldots,B$}

    \State
    \[
    r_t^{\splus}
    \gets
    \dfrac{w_t^{\splus}}{\lambda_t^{\splus}} +
    cA_t^*
    \left[
    A(z^{\splus}-y)-A_t(z_t^{\splus}-y_t)
    \right]
    -
    \frac{1}{\lambda_t^{\splus}}(z_t^{\splus}-y_t).
    \]
    \label{ln:PABIPP:residual}

\EndFor

%\Statex

\State $r^{\splus}\gets(r_1^{\splus},\ldots,r_B^{\splus})$ and $\delta^{\splus} \gets\sum_{t=1}^B
(\epsilon_t^{\splus}/\lambda_t^{\splus}).
$\label{ln:PABIPP:assemble}

\State \Return $(y^{\splus},\lambda^{\splus};z^{\splus},r^{\splus},\delta^{\splus})$

\end{algorithmic}
\end{algorithm}
\addtocounter{algorithm}{-1}
\renewcommand{\thealgorithm}{\arabic{algorithm}}
}

We next clarify two details of {\PABIPP}. First, we comment on the inexact solution of the block subproblems and the adaptive selection of the proximal stepsizes. The parallel \textbf{for} of lines~\ref{ln:PABIPP:for}--\ref{ln:PABIPP:endfor} performs $B$ independent blockwise procedures. For each $t\in\{1,\ldots,B\}$, we refer to the procedure associated with block $t$ as the $t$-th {\PABIPP} loop. Each iteration of the $t$-th {\PABIPP} consists of two main steps: computing an inexact solution to the $t$-th block subproblem in line~\ref{ln:PABIPP:solve} and testing the corresponding block sufficient-decrease condition in line~\ref{ln:PABIPP:test}. The subroutine for finding an approximate solution of \eqref{algo.PABIPP-subproblem} as in Definition \ref{def:rel.stationary} will be discussed in the second paragraph after Proposition~\ref{prop:PABIPP}. Once the inexact solution is obtained, we test inequality~\eqref{algo:PABIPP.adaptive.lambda} with the current trial stepsize $\lambda_t^+$. If this inequality holds, $\lambda_t^{\splus}$ is halved in line~\ref{ln:PABIPP:halve}, and the loop returns to line~\ref{ln:PABIPP:solve}, where the block subproblem is solved again with the new stepsize. Otherwise, we accept the current solution and stepsize, and the $t$-th {\PABIPP} loop terminates.

Second, {\PABIPP} generates two primal points with complementary roles. The point \(z^{\splus}\) is associated with the inexact stationary solutions of the block proximal subproblems. In contrast, \(y^{\splus}\) is constructed to satisfy the required decrease condition. Although both primal points $z^{\splus}$ and $y^{\splus}$ are passed to the outer algorithm, $z^{\splus}$ is used to update the Lagrange multiplier, while $y^{\splus}$ is used as the next primal iterate. The precise properties of these points and their roles in the outer algorithm are described right after Proposition~\ref{prop:PABIPP}.

Before stating the main result of this section, we define the following quantities 
\begin{align}\label{eq:def.zeta2andblock_norm}
\zeta:= 8\bar\tau^{-2}\|A\|_\dagger^2, \qquad \|A\|_\dagger := \sqrt{\sum_{t=1}^B \|A_t\|^2},
\end{align}
with submatrices $A_t$'s as in \eqref{eq:block_structure}.

\begin{proposition}\label{prop:PABIPP}
For a given quadruple $(y,q,\lam,c)
\in {\cal H} \times A(\mathbb{R}^n) 
\times \R^B_{++}\times \R_{++}$, the following statements about {\PABIPP} with input $(y,q,\lam, c)$ hold:

\begin{itemize}
\item[(a)] for every $t \in \{1,\ldots,B\}$,  the $t$-th {\PABIPP} loop stops  in at most $1+\lceil \log(1+2m_t\lam_t) \rceil$ iterations and the prox stepsize $\lam^{\splus}$
output by {\PABIPP} satisfies
\begin{equation}\label{eq:ABIPP:output.lambda.bound.2}
\lam_t^{\splus} \le \lam_t, \quad \max \left\{ \frac{1}{\lambda_t}, 4 m_t \right\} = 
 \max \left\{ \frac{1}{\lambda_t^{\splus}}, 4 m_t \right\}
  \quad \forall t=1,\ldots,B;
\end{equation}

\item[(b)] if $(y^{\splus},\lam^{\splus};z^{\splus}, r^{\splus},\delta^{\splus})={\PABIPP}(y,q,\lam, c)$, then  
\begin{align}
r^{\splus}  &\in \nabla f (z^{\splus})+ \partial_{\delta^{\splus}} h(z^{\splus}) +
 A^*[q+c(Az^{\splus}-b)], \label{prop.B-IPP.inclusion} \\
 \|r^{\splus}\|^2 + \delta^{\splus}  
&\le \left[\frac{1}{16\bar\tau^2}\left(\frac{130}{\min(\lambda^{\splus})}+1\right) + c \zeta\right]
  \Big[ {\cal L}_c(y;q) - {\cal L}_c(y^{\splus};q)\Big].\label{prop.B-IPP.inequality}
\end{align}
where $\bar \tau=1/[4(B-1)]$ is as in line \ref{ln:PABIPP:tau} of {\PABIPP}.
\end{itemize}
\end{proposition}

% We postpone the proof of Proposition \ref{prop:ABIPP.output} to the end of Subsection \ref{prop:ABIPP.output}, and now only make a few remarks about it.
% First, statement (a) guarantees that for every $t\in\{1,\ldots,B\}$, the $t$-th {\PABIPP} loop terminates after a finite, logarithmically bounded number of iterations.
% Moreover, relation \eqref{eq:ABIPP:output.lambda.bound.2} shows that the accepted prox stepsize $\lam_t^{\splus}$ is no larger than the initial stepsize $\lam_t$, while preventing the backtracking procedure from reducing it unnecessarily. \textcolor{red}{Improve: When $m_t=0$, relation \eqref{eq:ABIPP:output.lambda.bound.2} directly gives $\lambda_t^+=\lambda_t$.}

We defer the proof of Proposition~\ref{prop:PABIPP} to Section~\ref{sec:proof.prop}, and limit the discussion here to two comments concerning this proposition. First, statement~\textit{(a)} guarantees that for every $t\in\{1,\ldots,B\}$, the $t$-th {\PABIPP} loop terminates. More specifically, this loop performs at most $1+\lceil\log(1+2m_t\lam_t)\rceil$ iterations, and hence halves the prox stepsize at most $\lceil\log(1+2m_t\lam_t)\rceil$ times. Moreover, it follows from the proof of statement~\textit{(a)} that this loop performs a single iteration, and hence $\lam_t^{\splus}=\lam_t$, whenever $2m_t\lam_t\le 1$ (in particular, whenever $m_t=0$, i.e., $f_t$ is convex on ${\cal H}_t$). 
Second, the inclusion \eqref{prop.B-IPP.inclusion} shows that $(r^{\splus},\delta^{\splus})$ is a residual pair for the primal point $z^{\splus}$.
Third, inequality \eqref{prop.B-IPP.inequality} bounds the magnitude of this residual pair in terms of the variation ${\cal L}_c(y;q) - {\cal L}_c(y^{\splus};q)$ of the augmented Lagrangian. 
This variation will play the role of a potential decrease in the analysis of the main algorithm. It is worth noting that the residual pair is associated with the primal point $z^{\splus}$, whereas the Lagrangian decrease is evaluated at $y^{\splus}=y+\bar\tau(z^{\splus}-y).$

We conclude this subsection with a few remarks on how to compute, for any given $t\in\{1,\ldots,B\}$, an $(1/32;y_t)$-relative stationary solution of \eqref{algo.PABIPP-subproblem}, as required in line \ref{ln:PABIPP:solve} of {\PABIPP}.
First, if an exact solution $z_t^{\splus}$ of \eqref{algo.PABIPP-subproblem} can be computed in closed form, then $(z_t^{\splus},w_t^{\splus},\epsilon_t^{\splus})=(z_t^{\splus},0,0)$ is a $(1/32;y_t)$-relative stationary solution of that subproblem.

Second, although the complexity analysis of {\PABIPP} does not
require the gradients $\nabla f_t$ to be Lipschitz continuous,
we impose this additional assumption to discuss a possible
implementation of the subproblem solver. Specifically, suppose
that $\nabla f_t$ is $\widetilde L_t$-Lipschitz continuous with
respect to the $t$-th block $x_t$. Then the differentiable part
of the objective function in \eqref{algo.PABIPP-subproblem}
has a Lipschitz continuous gradient with constant bounded by $1+\lambda_t^{\splus}\bigl(\widetilde L_t+c\|A_t\|^2\bigr).$
Under this additional assumption, the ADAP-FISTA procedure
described in \cite[Appendix B]{MaiaEtAl2026}, initialized with
$
(z^0,M_0,\mu_0,\sigma)
=
(y_t,\lam_t^{\splus} c\|A_t\|^2,1/2,1/1024),
$
can be used to find the required relative stationary solution
of \eqref{algo.PABIPP-subproblem}.
Whenever the procedure terminates successfully, its output satisfies this relative stationarity criterion.
The possibility of failure is discussed next.

Third, ADAP-FISTA may
fail to find the required relative stationary solution if
$\lambda_t^{\splus}>1/(2m_t)$ in line~\ref{ln:PABIPP:solve}
of {\PABIPP}. Indeed, in this case, the differentiable part
of the subproblem objective is not necessarily
$(1/2)$-strongly convex
(see \cite[Lemma 6.4]{MaiaEtAl2026}), and hence the procedure
may fail when invoked with $\mu_0=1/2$
(see \cite[Proposition B.1(c)]{MaiaEtAl2026}).
Such a failure is treated as a signal to reduce the current
proximal stepsize: $\lambda_t^{\splus}$ is halved regardless of
whether \eqref{algo:PABIPP.adaptive.lambda} is satisfied.
In either case, whether successful or not, ADAP-FISTA
terminates within
\[
{\cal O}\!\left(
\left[1+\lambda_t^{\splus}
\bigl(\widetilde L_t+c\|A_t\|^2\bigr)\right]^{1/2}
\right)
\]
iterations. Although this bound depends on $\widetilde L_t$, the procedure does not require prior knowledge of this constant.

% If an exact solution $u_t^{\mathrm{a}}$ of \eqref{algo.PABIPP-subproblem} can be computed in closed form, then $(u_t^{\mathrm{a}},r_t^{\mathrm{a}},\varepsilon_t^{\mathrm{a}})=(u_t^{\mathrm{a}},0,0)$ is a $(1/32;y_t)$-relative stationary solution of \eqref{algo.PABIPP-subproblem}.
% More generally, an inexact composite optimization routine may be used to generate the required relative stationary triple. 
% For example, whenever the \textcolor{blue}{differentiable} part of \eqref{algo.PABIPP-subproblem} satisfies the assumptions required by ADAP-FISTA, see \cite[Appendix B]{MaiaEtAl2026}, the latter can be applied with relative accuracy parameter $\sigma=1/\sqrt{32}$ and the current trial stepsize $\lambda_t^+$. 
% If $\nabla_t f$ is $\widetilde L_t$-Lipschitz continuous with respect to the $t$-th block, then the gradient of the \textcolor{blue}{differentiable} part of \eqref{algo.PABIPP-subproblem} has Lipschitz constant bounded by
% $
% 1+\lambda_t^+\bigl(\widetilde L_t+c\|A_t\|^2\bigr).
% $
% As in the discussion following \cite[Proposition 3.1]{MaiaEtAl2026}, the adaptive nature of ADAP-FISTA means that $\widetilde L_t$ need not be known in advance. Finally, \cite[Proposition 4.1]{MaiaEtAl2026} shows that the number of iterations performed by ADAP-FISTA is bounded by ${\cal O} ( [\tilde L_t + c \|A_t\|^2 ]^{1/2})$.

%\textcolor{red}{Comments about Proposition 4.1}

\subsection{The {\PAADMM}}
\label{sub:PAADMM}

This subsection presents {\PAADMM} and states its main complexity result. We first describe the method and discuss its main properties, and then establish a bound on the total number of iterations to obtain a near-stationary solution of \eqref{initial.problem}--\eqref{eq:block_structure} in the sense of \eqref{eq:rho.eta.solution}.

\begin{algorithm}[H]
\caption{{\PAADMM}}
\begin{algorithmic}[1]
\Statex \hskip-1.8em
\textbf{Universal Input:} $(\rho,\eta)\in\R^2_{++}$, $\alpha>0$, $C\ge\rho$
\Require $x^0\in{\cal H}$ and $\gamma^0\in\R^B_{++}$
\Ensure $(x,p,v,\varepsilon,c)$ that satisfies the conclusions of Theorem~\ref{thm:main}.

\Statex

\State $p^0\gets 0\in\R^l$,
$c_0\gets 1/[1+\|Ax^0-b\|]$
\label{ln:PAADMM:init}
\Comment{initialization}

\For{$s\gets1,2,\ldots$}

    \State $T_0\gets0$, $R_0\gets1$, $k\gets0$, $c\gets c_{s-1}$, $(y^0,q^0,\lambda^0)\gets
    (x^{s-1},p^{s-1},\gamma^{s-1})$.
    \label{ln:PAADMM:stage.init}
    \Comment{stage initialization}

    \For{$i\gets1,2,\ldots$}

        \State $(y^i,\lambda^i;z^i,r^i,\delta_i)
        \gets
        \PABIPP(y^{i-1},q^{i-1},\lambda^{i-1},c)$ \label{ln:PAADMM:call}

        \State $\mathsf{q}^{\text{cand}} \gets q^{i-1} + c(Az^i-b)$ \Comment{candidate multiplier}
        
        \If{$\|r^i\|^2+\delta_i\le\rho^2$}
        \label{ln:PAADMM:stat.test}
        \Comment{stationarity residual bound check}

            \State $k\gets k+1$, \ $\cali_k\gets i$
            \label{ln:PAADMM:count.stop}

            \State $q^i\gets \mathsf{q}^{\text{cand}}$

            \State $(x^s,p^s,v^s,\varepsilon_s,\gamma^s)
            \gets
            ( z^i,q^i,r^i,\delta_i,\lambda^i)$
            and \textbf{go to line~\ref{ln:PAADMM:double.c}}\label{ln:PAADMM:stage.output}

        \EndIf

        \State $T_i
        \gets
        \Lc(y^{i-1};q^{i-1})
        -
        \Lc(y^i;q^{i-1})
        +
        T_{i-1}$\label{ln:PAADMM:T}
        
        \State $R_i = \dfrac{1}{
    \max
    \left\{
    \displaystyle \max_{1\leq j\leq i-1} 
    \{
    \|q^{j} - q^{j-1}\|^2
    \},
    \|\mathsf{q}^{\text{cand}}-q^{i-1}\|^2
    \right\}
    }
    \left[\displaystyle \sum_{j=1}^{i-1}\|q^{j}-q^{j-1}\|^2 + \|\mathsf{q}^{\text{cand}}-q^{i-1}\|^2 \right]$ \label{ln:PAADMM:R}

        \Statex

        \Statex
        
        \If{$\|r^i\|^2+\delta_i\le C^2$
        and
        $\dfrac{\frac{1}{i}\sum_{j=1}^{i}(\|r^j\|^2+\delta_j)}{\alpha R_i}
        \ge
        \dfrac{T_i}{i}$}
        \label{ln:PAADMM:LMT}
        \Comment{Lagrange multiplier test (LMT)}

            \State $k\gets k+1$, \ $\cali_k\gets i$
            \label{ln:PAADMM:count.LMT}

            \State $q^i\gets \mathsf{q}^{\text{cand}}$
            \label{ln:PAADMM:mult.LMT}

        \Else

            \State $q^i\gets q^{i-1}$
            \label{ln:PAADMM:mult.keep}

        \EndIf

    \EndFor

    \State $c_s\gets2c_{s-1}$
    \label{ln:PAADMM:double.c}

    \If{$\|Ax^s-b\|\le\eta$}
    \label{ln:PAADMM:feas.test}
    \Comment{feasibility bound check}

        \State \Return $(x,p,v,\varepsilon,c)
        =(x^s,p^s,v^s,\varepsilon_s,c_s)$

    \EndIf \label{ln:PAADMM:stage.end}

\EndFor

\end{algorithmic}
\end{algorithm}

A few comments are in order about {\PAADMM}. First, {\PAADMM} consists of two main loops. The outer loop is indexed by $s$ and will be referred to as a \textit{stage}. During \textit{stage} $s$, the penalty parameter is kept fixed and equal to $c_{s-1}$. The stage starts in line~\ref{ln:PAADMM:stage.init} and ends in line~\ref{ln:PAADMM:stage.end}.

Second, within each stage, {\PAADMM} performs a sequence of \textit{inner iterations}, indexed by $i$. Each inner iteration performs one call to {\PABIPP} in line~\ref{ln:PAADMM:call}. The index $k$, initialized in line~\ref{ln:PAADMM:stage.init} and incremented in lines~\ref{ln:PAADMM:count.stop} and~\ref{ln:PAADMM:count.LMT}, counts the number of Lagrange multiplier updates performed within the current \textit{stage}. The index $\cali_k$, computed either in line~\ref{ln:PAADMM:count.stop} or line~\ref{ln:PAADMM:count.LMT}, denotes the inner iteration at which the $k$-th Lagrange multiplier update occurs. It follows from Theorem~\ref{thm:main}(b) below that the total number of inner iterations performed by {\PAADMM} is finite, and hence $\cali_k$ is well-defined. If the inequality in line~\ref{ln:PAADMM:stat.test} is satisfied, the inner loop performs its final Lagrange multiplier update and terminates in line~\ref{ln:PAADMM:stage.output}, returning the quintuple $(x^s,p^s,v^s,\varepsilon_s,\gamma^s)$. Otherwise, depending on the LMT in line~\ref{ln:PAADMM:LMT}, the method either updates the Lagrange multiplier in line~\ref{ln:PAADMM:mult.LMT} or leaves it unchanged in line~\ref{ln:PAADMM:mult.keep}. In both cases, it then proceeds to the next inner iteration.

Third, once the output $(x^s,p^s,v^s,\varepsilon_s,\gamma^s)$ from the current \textit{stage} has been obtained, the penalty parameter is doubled in line~\ref{ln:PAADMM:double.c}, and feasibility is tested in line~\ref{ln:PAADMM:feas.test}. If this test is satisfied, the algorithm terminates; otherwise, a new \textit{stage} is initiated.

Finally, {\PAADMM} is fully adaptive, i.e., it does not require prior knowledge of the instance parameters associated with problem~\eqref{initial.problem}--\eqref{eq:block_structure}. Indeed, {\PABIPP} adaptively selects the proximal stepsize (see Subection~\ref{sub:PABIPP} and Section~\ref{sec:proof.prop}), while the remaining components of {\PAADMM} do not require knowledge of either Lipschitz constants or weak convexity parameters.

Before stating the main results of this section, namely Theorem~\ref{thm:main} and Corollary~\ref{cor:complexity}, we define the quantities
\begin{equation}\label{eq:tildezeta.Upsilon}
    \tilde \zeta =\frac{1}{16\bar\tau^2}\left[ 130\max_{1\leq t\leq B}\left\{\frac{1}{\gamma_t^0} + 4m_t\right\}+1\right], \qquad \Upsilon(C):=
\frac{2D_hM_h+(2D_h+1)(C+C^2+\widehat\nabla_f)}
{\bar d\,\nu_A^+}  + \textcolor{red}{1}, 
\end{equation}
and
\begin{align}\label{eq:def.Gamma.complexity}
\Gamma(x^0;C,\alpha)&:=2(\phi_{\mathrm{sup}}-\phi_{\mathrm{inf}}) +(1+\|Ax^0-b\|)\Upsilon(C)^2\left[ 36 + \frac{16\tilde \zeta }{\alpha}\left(\frac{4B(1-\bar\tau)^2}{\bar\tau}+1\right)\right] \nonumber\\
 &\qquad\qquad + \frac{16\zeta \Upsilon(C)^2}{\alpha}\left(\frac{4B(1-\bar\tau)^2}{\bar\tau}+1\right),
\end{align}
where $\bar \tau$ is as in line \ref{ln:PABIPP:tau} of {\PABIPP}, $x^0$ is the input of {\PAADMM}, $M_h$ and $\bar d$ are as in assumptions (A1) and (A3), respectively, $(D_h,\widehat\nabla_f)$ is as in \eqref{def:damH}, and $\nu_A^+$ denotes the smallest positive singular value of the nonzero linear operator $A$.

\begin{theorem}
\label{thm:main}
For a given universal input $(\rho,\eta,\alpha,C)$, with $(\rho,\eta)\in\R^2_{++}$, $\alpha>0$, and $C\ge\rho$, and a given input $(x^0,\gamma^0)\in{\cal H}\times \R_{++}^B$, the following statements about {\PAADMM} with these inputs hold:
\begin{enumerate}
    \item[(a)] for every stage $s$ generated by {\PAADMM}, the quadruple $(x^{s},p^{s},v^s,\varepsilon_s)$ obtained at the end of the $s$-th stage is a $\rho$-stationary solution of \eqref{initial.problem}--\eqref{eq:block_structure}, i.e., it satisfies the first two conditions in \eqref{eq:rho.eta.solution}, and the following bounds hold
\begin{equation}
\label{eq:prop:input.bounds}
\|p^{s}\| \le \Upsilon(C) \quad \text{and} \quad  c_{s}\|Ax^{s} -b \|
\le 4\Upsilon(C),
\end{equation}
where $\Upsilon(C)$ is as in~\eqref{eq:tildezeta.Upsilon};
\item[(b)] for every stage $s$ generated by {\PAADMM}, the number of inner iterations performed within the $s$-th stage is bounded by  
\begin{align}\label{eq:prop.inner.loop.totaliterations}
 \left(\frac{\tilde\zeta+c_{s-1}\zeta}{\rho^2}\right)\Gamma(x^0; C,\alpha) + 1,
\end{align}
where  $\zeta$, $\tilde\zeta$ and $\Gamma(x^0; C,\alpha)$ are as in \eqref{eq:def.zeta2andblock_norm}, \eqref{eq:tildezeta.Upsilon} and \eqref{eq:def.Gamma.complexity}, respectively.

\item[(c)] {\PAADMM} terminates after at most $\log[1+4\Upsilon(C)/(c_0\eta)]+1$ stages with a $(\rho,\eta)$-stationary solution of \eqref{initial.problem}--\eqref{eq:block_structure}.
\end{enumerate}
\end{theorem}

We defer the proof of Theorem~\ref{thm:main} to the end of technical Section~\ref{sec:proof.thm}, and now just make some comments about it. First, statement~(a) shows that each stage of {\PAADMM} ends with a $\rho$-stationary solution of \eqref{initial.problem}--\eqref{eq:block_structure} that also satisfies the two important inequalities in \eqref{eq:prop:input.bounds}. The first inequality shows that the multipliers remain bounded by $\Upsilon(C)$ in every stage. The second inequality shows that $\|Ax^s-b\|\le 4\Upsilon(C)/c_{s}$, so doubling the penalty parameter at each stage (see line~\ref{ln:PAADMM:double.c} of {\PAADMM}) will eventually make $x^s$ $\eta$-feasible.

Second, since $\bar\tau=1/[4(B-1)]$ by line~\ref{ln:PABIPP:tau} of {\PABIPP}, it follows from \eqref{eq:def.zeta2andblock_norm} and \eqref{eq:tildezeta.Upsilon} that $\zeta=\Theta(B^2)$ and $\tilde\zeta=\Theta(B^2)$. Moreover, it follows from \eqref{eq:def.Gamma.complexity} that $\Gamma(x^0;C,\alpha)={\cal O}(1)$ under the mild assumptions that $\|Ax^0-b\|={\cal O}(1)$, $\alpha=\Omega(B^4)$, and $C={\cal O}(1)$. Hence, statement~(b) shows that the $s$-th stage performs ${\cal O}(B^2(1+c_{s-1})\rho^{-2})$ inner iterations.

Finally, statement~(c) shows that the penalty parameter upon termination is 
${\cal O}(\eta^{-1})$, 
and hence that 
{\PAADMM} performs 
at most ${\cal O}(\log(2+\eta^{-1}))$ penalty updates to obtain a $(\rho,\eta)$-stationary solution of 
\eqref{initial.problem}--\eqref{eq:block_structure}.  

As a consequence of Theorem~\ref{thm:main}, we now establish the main result of this paper, namely, the iteration complexity of {\PAADMM}, i.e., the total number of inner iterations it needs to obtain a $(\rho,\eta)$-stationary solution of ~\eqref{initial.problem}--\eqref{eq:block_structure}.

\begin{cor}\label{cor:complexity}
({\PAADMM} Iteration Complexity) The total number of inner iterations performed by {\PAADMM} is bounded by
\begin{equation*}%\label{com.AADMM}
\left[
1+\frac{\tilde\zeta\Gamma(x^0;C,\alpha)}{\rho^2}
\right]
\log\left(2+\frac{8\Upsilon(C)}{c_0\eta}\right)
+
\frac{\zeta\Gamma(x^0;C,\alpha)}{\rho^2}
\left[
c_0+\frac{8\Upsilon(C)}{\eta}
\right]
\end{equation*} 
where $\zeta$, $\tilde\zeta$, $\Upsilon(C)$   and $ \Gamma(x^0;C,\alpha)$ are as in~\eqref{eq:def.zeta2andblock_norm}, \eqref{eq:tildezeta.Upsilon}, and \eqref{eq:def.Gamma.complexity}, respectively, and $c_0$ is as in line \ref{ln:PAADMM:init} of {\PAADMM}.    
\end{cor}

\begin{proof}
Let $\bar s$ denote the total number of stages generated by {\PAADMM}. By Theorem~\ref{thm:main}(c), we have
\begin{equation}\label{def:bar:s}
\bar s
\le
1+\log\left(1+\frac{4\Upsilon(C)}{c_0\eta}\right)
=
\log\left(2+\frac{8\Upsilon(C)}{c_0\eta}\right).
\end{equation}
Hence, summing the bound~\eqref{eq:prop.inner.loop.totaliterations} of Theorem~\ref{thm:main}(b) over $s=1,\ldots,\bar s$ and using~\eqref{def:bar:s}, we conclude that the total number of inner iterations performed by {\PAADMM} is bounded by
\begin{align*}
\sum_{s=1}^{\bar s}
\left[
\left(
\frac{\tilde\zeta+c_{s-1}\zeta}{\rho^2}
\right)
\Gamma(x^0;C,\alpha)
+1
\right]
&=
\left[
1+\frac{\tilde\zeta\Gamma(x^0;C,\alpha)}{\rho^2}
\right]\bar s
+
\frac{\zeta\Gamma(x^0;C,\alpha)}{\rho^2}
\sum_{s=1}^{\bar s}c_{s-1}\\
&\leq \left[
1+\frac{\tilde\zeta\Gamma(x^0;C,\alpha)}{\rho^2}
\right]
\log\left(2+\frac{8\Upsilon(C)}{c_0\eta}\right)
+
\frac{\zeta\Gamma(x^0;C,\alpha)}{\rho^2}\sum_{s=1}^{\bar s}c_{s-1}.
\end{align*}

Using that $c_{s-1}=2^{s-1}c_0$ and \eqref{def:bar:s}, we obtain that
\begin{align*}
\sum_{s=1}^{\bar s}c_{s-1}
=c_0\sum_{s=1}^{\bar s}2^{s-1}=c_0(2^{\bar s}-1)\le c_0\left(2+\frac{8\Upsilon(C)}{c_0\eta}-1\right)= 
c_0+\frac{8\Upsilon(C)}{\eta}.
\end{align*}
Hence, combining the two previous inequalities completes the proof.
\end{proof}

We now write the dependence of the bound in Corollary~\ref{cor:complexity} in terms of the tolerances $(\rho,\eta)$ and on the number of blocks $B$. Recall from the second paragraph after Theorem~\ref{thm:main} that if we assume that $\|Ax^0-b\|={\cal O}(1)$, $C=\Theta(1)$, and $\alpha\ge B^4$, then $\zeta={\cal O}(B^2)$, $\tilde\zeta={\cal O}(B^2)$, $\Upsilon(C)={\cal O}(1)$, and $\Gamma(x^0;C,\alpha)={\cal O}(1)$. Hence, the bound in Corollary~\ref{cor:complexity} becomes
\[
{\cal O}\left(\frac{B^2}{\rho^2}\log\Big(2+\frac{1}{\eta}\Big)+\frac{B^2}{\rho^2\eta}\right)={\cal O}\left(\frac{B^2}{\rho^2\eta}\right).
\]
In particular, with $\epsilon:=\min\{\rho,\eta\}$, {\PAADMM} finds a $(\rho,\eta)$-stationary solution in ${\cal O}(B^2\epsilon^{-3})$ inner iterations.

\section{Proof of Proposition \ref{prop:PABIPP}}\label{sec:proof.prop}

This section establishes the properties of {\PABIPP} stated in
Proposition~\ref{prop:PABIPP}.

\begin{lemma}\label{lem:jacobi.decrease}
Let $y\in{\cal H}$, $z^{\splus}=(z_1^{\splus},\ldots,z_B^{\splus})\in{\cal H}$, $(q,c)\in A(\R^n)\times\R_{++}$, and $\tau\in[0,1]$ be given, and define
\begin{equation}\label{def:z_tau}
y(\tau):=y + \tau(z^{\splus}-y). 
\end{equation} 
Then,
\begin{align}\label{eq:lemma:jacobi-gauss}
\Lc(y;q) - \Lc\left(y(\tau);q\right) &\geq \sum_{t=1}^{B}\Big[\Lc(y;q) - \Lc(y_{<t},y_{t}(\tau),y_{>t};q)\Big] 
- \frac{c(B-1)\tau^2}{2}\sum_{t=1}^{B}\|A_t(z_t^{\splus}-y_t)\|^2.
\end{align}
\end{lemma}

\begin{proof}
Using Lemma \ref{lem:AL.block.identity} with $(x,y)=(y,y(\tau))$ (note that $y(\tau)\in{\cal H}$ because ${\cal H}$ is convex and $\tau\in[0,1]$) and the identity $y_t(\tau)-y_t=\tau(z_t^{\splus}-y_t)$, we have that for every $t\in\{1,\ldots,B\}$,
\begin{align*}
    \Lc(y_{<t}(\tau),y_t,y_{> t};q) - &\Lc(y_{< t}(\tau),y_{t}(\tau),y_{>t};q) \\
    &= \Lc(y;q) - \Lc(y_{<t},y_{t}(\tau),y_{>t};q) - c\tau^2\sum_{i=1}^{t-1}\inner{A_t(z_t^{\splus}-y_t)}{A_i(z_i^{\splus}-y_i)}.
\end{align*}
Summing the above equation over $t=1,\ldots,B$, we obtain
\begin{align}\label{eq:jacobi-gauss.proof}
\Lc(y;q) - \Lc(y(\tau);q) &\geq \sum_{t=1}^{B}\Big[\Lc(y;q) - \Lc(y_{<t},y_{t}(\tau),y_{>t};q)\Big] 
- c\tau^2\sum_{t=1}^B\sum_{i=1}^{t-1}\inner{A_t(z_t^{\splus}-y_t)}{A_i(z_i^{\splus}-y_i)}.
\end{align}
Now, for every $t\in\{1,\ldots,B\}$ and every $i=1,\ldots,t-1$, the Cauchy--Schwarz and AM--GM inequalities imply that
\[
\inner{A_t(z_t^{\splus}-y_t)}{A_i(z_i^{\splus}-y_i)} \leq \|A_t(z_t^{\splus}-y_t)\| \|A_i(z_i^{\splus}-y_i)\| \leq\frac{\|A_t(z_t^{\splus}-y_t)\|^2+\|A_i(z_i^{\splus}-y_i)\|^2}{2}.
\]
Hence, substituting the two previous inequalities into \eqref{eq:jacobi-gauss.proof}, we conclude that ~\eqref{eq:lemma:jacobi-gauss} holds.
\end{proof}

We next prove some properties of the proximal block subproblem when
the prox stepsize is sufficiently small.

\begin{lemma}\label{lem:Psi.convexity}
For a given $t\in\{1,\ldots,B\}$, let $\Psi_t(\cdot)$ denote the differentiable part of the objective function of the $t$-th block subproblem \eqref{algo.PABIPP-subproblem}, i.e.,
\begin{equation}\label{def:Psi:t}
\Psi_t(\cdot):=\lam_t^{\splus} {\cal L}_{c}^{(\mathrm{s})}(y_{<t},\cdot,y_{>t};q)+\frac{1}{2}\|\cdot-y_{t}\|^{2}
\end{equation}
where ${\cal L}_{c}^{(\mathrm{s})}(\cdot\,;q)$ is as in \eqref{def:smooth:ALM}. Moreover, let $B_t : = (1-\lam_t^{\splus} m_t)I+\lam_t^{\splus} cA_t^*A_t$ and assume that $1/\lam_t^{\splus} \in [2m_t,\infty)$. Then, the following statements hold:

\begin{itemize}
    \item[(a)] the matrix $B_t$ satisfies
    \begin{equation}\label{eq:B_t.positive}
\langle\;\cdot\;, \, B_t(\cdot) \, \rangle 
\ge 
\lam_t^{\splus} c \|A_t (\cdot)\|^2 +
 \frac12 \|\cdot\|^2,
    \end{equation}
%     \begin{equation}\label{eq:B_t.positive}
% \langle\;\cdot\;, \, B_t(\cdot) \, \rangle -\frac12 \|\cdot\|^2
% = 
% \left( \frac12 - \lam_t m_t \right) \|\cdot\|^2 + \lam_t^{\splus} c  \|A_t (\cdot)\|^2,
%     \end{equation}
and hence $\|\cdot\|_{B_t} = \sqrt{\inner{\cdot}{B_t( \cdot)}}$ defines a norm; 
    \item[(b)] both functions
    \begin{equation}\label{eq.Psi_t.strg.cvx}
        \Psi_t (\cdot) - \frac12 \|\cdot\|^2_{B_t}, \quad \Psi_t (\cdot) - \frac14 \|\cdot\|^2
    \end{equation}
    are convex on ${\cal H}_t$.
\end{itemize}
\end{lemma}
\begin{proof}
(a) Inequality \eqref{eq:B_t.positive} follows immediately from the fact that $1/\lam_t^{\splus} \in [2m_t,\infty)$. Moreover, it follows from \eqref{eq:B_t.positive} that $B_t$ is a positive definite matrix, and hence $\|\cdot\|_{B_t} = \sqrt{\inner{\cdot}{B_t( \cdot)}}$ defines a norm.  

(b) Using assumption (A2), the definition of the norm $\|\cdot\|_{B_t}$, the definitions of ${\cal L}_{c}^{(\textrm{s})}(\cdot\,;\cdot)$ and  $\Psi_t(\cdot)$ in \eqref{def:smooth:ALM}
 and \eqref{def:Psi:t}, respectively, we see that 
 \[
 \Psi_t(\cdot) - \frac12 \|\cdot\|^2_{B_t} - 
\lam_t^{\splus} \left[ f(y_{<t},\cdot,y_{>t})+ \frac{ m_t}2  \|\cdot\|^2 \right]
 \]
is affine. 
Using the fact that $f(y_{<t},\cdot,y_{>t})$ is $m_t$-weakly convex on ${\cal H}_t$ and the previous observation, we conclude that $\Psi_t(\cdot) - \frac12 \|\cdot\|^2_{B_t}$ is convex on ${\cal H}_t$. Next, using the previous fact, the condition that $1/\lam_t^{\splus} \in [2m_t,\infty)$, and the identity
 \[
 \Psi_t(\cdot) - \frac14 \|\cdot\|^2 =  \left(\Psi_t(\cdot) - \frac12 \|\cdot\|^2_{B_t}\right) + \frac{1}{2}\left( \|\cdot\|^2_{B_t} -  \frac{1}{2}\|\cdot\|^2\right),
\]
we conclude that $\Psi_t(\cdot) - \frac14 \|\cdot\|^2$ is the sum of two convex functions, and hence it is convex on ${\cal H}_t$.
\end{proof}

% \begin{framed}
% \textcolor{red}{FOR AUTHORS, ONLY:} 
% \begin{align*}
% &\Psi_t(\cdot) - \lam_t^+ \left[ f(y_{<t},\cdot,y_{>t})+ \frac{ m_t}2  \|\cdot\|^2 \right] \\
% &=\lam_t^{\splus} \left[ \inner{A(y_{<t},\cdot,y_{>t})-b}{q} +\frac{c}{2}\|A(y_{<t},\cdot,y_{>t})-b\|^2\right]+\frac{1}{2}\|\cdot-y_{t}\|^{2} - \lam_t^+ \frac{ m_t}2  \|\cdot\|^2\\
% &=\lam_t^{\splus} \left[ \inner{A(y_{<t},\cdot,y_{>t})-b}{q} +\frac{c}{2}\|A(y_{<t},\cdot,y_{>t})-b\|^2\right]+\frac{1}{2}\|\cdot-y_{t}\|^{2}- \lam_t^+ \frac{ m_t}2  \|\cdot\|^2\\
% &=\lam_t^{\splus} \left[\frac{c}{2}\|A_t(\cdot)\|^2\right]+\frac{1-\lam_t^+ m_t}{2}\|\cdot\|^{2} + a(\cdot)\\
% & = \frac{1}{2}\|\cdot \|_{B_t}^2  + a(\cdot)
% \end{align*}
% \end{framed}

% \begin{framed}
% [RENATO] If $1/\lam_t^{\splus} \in [2m_t,\infty)$, then the norm $\|\cdot\|_{B_t}$ is well defined and the function $\Psi_t(\cdot) - \frac12 \|\cdot\|^2_{B_t}$ is convex. As a consequence, $\Psi_t(\cdot) - \frac12 \|\cdot\|^2$ is convex.

% 1) the function
% \[
% \Psi_t (\cdot) - \frac12 \inner{\cdot}{B_t \, \cdot}
% \]
%  is convex
 
% 2) if $\lam_t^+m_t \le 1/2$ then
%  \[
% \langle\;\cdot\;, \, B_t(\cdot) \, \rangle 
% \ge 
% \lam_t^{\splus} c \|A_t (\cdot)\|^2 +
%  \frac12 \|\cdot\|^2.
%  \]
% %  and hence 
% % \[
% % B_t \succeq (1/2) I + \lam_t^+ c A_t^*A_t
% % \]
% and hence $\|\cdot\|_{B_t} = \sqrt{\inner{\cdot}{B_t( \cdot)}}$ is a norm; 

% 3) both functions
% \[
% \Psi_t (\cdot) - \frac12 \|\cdot\|^2_{B_t}, \quad \Psi_t (\cdot) - \frac12 \|\cdot\|^2
% \]
% are convex;
% \end{framed}

Lemma~\ref{lem:Psi.convexity} gives the strong convexity required to
turn the inexact stationarity condition from the block solver
into a decrease estimate. More precisely, line~\ref{ln:PABIPP:solve} of {\PABIPP}, together with Definition~\ref{def:rel.stationary}, generates a triple
$(z_t^{\splus},w_t^{\splus},\epsilon_t^{\splus})$ satisfying
\begin{gather}\label{eq:ref:lemma:A1}
\begin{gathered}
w_t^{\splus}\in \nabla\left[\lam_t^{\splus} {\cal L}_{c}^{(\textrm{s})}(y_{<t},\cdot ,y_{>t}; q)+\frac{1}{2}\|\cdot -y_{t}\|^2\right](z_t^{\splus})+\partial_{\epsilon_t^{\splus}}(\lam_t^{\splus} h_t)(z_t^{\splus}), \\
4\|w_t^{\splus}\|^2 + 2 \epsilon_t^{\splus} \leq \frac{1}{16}\| z_t^{\splus} -y_t\|^2,
\end{gathered}
\end{gather}
where ${\cal L}_{c}^{(\textrm{s})}(\, \cdot \, ; q)$ is  as in \eqref{def:smooth:ALM}. Then, Lemma~\ref{lem:block.decrease} combines \eqref{eq:ref:lemma:A1} with Lemma~\ref{lem:Psi.convexity} to yield a decrease at the block candidate $z_t^{\splus}$. Throughout the analysis, we will also repeatedly use the functions $\Psi_t(\cdot)$, $t=1,\ldots,B$, defined in \eqref{def:Psi:t}.

\begin{lemma}\label{lem:block.decrease}
Let $t \in \{1,\ldots,B\}$. If the prox stepsize $\lam_t^{\splus}$ satisfies $1/\lam_t^{\splus} \in [2m_t,\infty)$, then
    \begin{align}\label{eq:lemma.block.decrease.hatz}
        {\cal L}_c(y; q)-{\cal L}_c(y_{<t}, z_t^{\splus}  ,y_{>t}; q)  &\ge 
 \frac{1}{2\lam_t^{\splus}} \|z_t^{\splus} - y_t\|^2 + \frac{c}{4} \|A_t (z_t^{\splus} - y_t)\|^2.
    \end{align}
\end{lemma}
\begin{proof}
We first claim that the quadruple $(z_t^{\splus},w_t^{\splus}, \epsilon_t^{\splus},\lam_t^{\splus})$ satisfies
\begin{align} \label{aux:inclusion}
w_t^{\splus}\in  \partial_{\epsilon_t^{\splus}}
         \left[ \Psi_t + \lam_t^{\splus} h_t\right](z_t^{\splus}),
\end{align}
where $\Psi_t(\cdot)$ is as in \eqref{def:Psi:t}.
Indeed, it follows from  \eqref{def:Psi:t} and the inclusion in \eqref{eq:ref:lemma:A1} that the inclusion $w_t^{\splus}\in  \nabla \Psi_t(z_t^{\splus}) + \partial_{\epsilon_t^{\splus}}(\lam_t^{\splus} h_t) (z_t^{\splus} )$ holds.
Moreover, by noting that $\Psi_t(\cdot)$ and $\lam_t^{\splus} h_t(\cdot )$ are convex functions (see Lemma \ref{lem:Psi.convexity}(b) and assumption (A1)) on ${\cal H}_t$, and using \cite[Thm. 3.1.1 of Ch.\ XI]{lemarechal1993}, we conclude that $\nabla \Psi_t(z_t^{\splus}) + \partial_{\epsilon_t^{\splus}}(\lam_t^{\splus} h_t) (z_t^{\splus} )\subset \partial_{\epsilon_t^{\splus}}\left[ \Psi_t + \lam_t^{\splus} h_t\right](z_t^{\splus})$. 
The claim now follows from the two previous inclusions.
% Hence, the quadruple $(u_t^{\mathrm{a}},r_t^{\mathrm{a}}, \varepsilon_t^{\mathrm{a}},\lam_t^{\splus})$ satisfy
% \begin{align*} r_t^{\mathrm{a}}\in  \partial_{\varepsilon_t^{\mathrm{a}}}
%          \left[ \Psi_t + \lam_t^{\splus} h_t\right](u_t^{\mathrm{a}}),
% \end{align*}

Thus, as the first function in \eqref{eq.Psi_t.strg.cvx} is convex and \eqref{aux:inclusion} holds, we can apply
%Now, using that $\Psi_t(\cdot) - \frac12 \|\cdot\|^2_{B_t}$ is  convex due to Lemma~\ref{lem:Psi.strgl.cvx}(b),
%we may therefore apply 
Lemma \ref{lem:eps.subgrad}  with $\psi=\Psi_t+\lambda_t^{\splus}h_t$, $(\xi,\varsigma, Q)=(1, 1, B_t)$, $(u,y,v) = (y_{t},z_t^{\splus},w_t^{\splus})$, and $\eta=\epsilon_t^{\splus}$, to obtain that
\begin{align*}
    &\lam_t^{\splus} {\cal L}_c(y_{<t},y_{t}, y_{>t}; q)-\left[\lam_t^{\splus} {\cal L}_c(y_{<t},z_t^{\splus}, y_{>t}; q) +\frac{1}{2}\|z_t^{\splus} - y_t\|^2 \right] \\
    &= \left[\Psi_t(y_t)+\lambda_t^{\splus} h_t(y_t)\right] - \left[ \Psi_t(z_t^{\splus})+\lambda_t^{\splus} h_t(z_t^{\splus} ) \right]  \\
    &\ge  \frac{1}{4}\|z_t^{\splus} - y_t\|_{B_t}^2 -2\epsilon_t^{\splus}+\langle w_t^{\splus}, y_t-z_t^{\splus}\rangle \overset{\eqref{eq:B_t.positive}}\ge \frac{1}{8}\|z_t^{\splus}-y_t\|^2 + \frac{\lam_t^{\splus} c}{4}\|A_t(z_t^{\splus}-y_t)\|^2 -2\epsilon_t^{\splus} +\langle w_t^{\splus}, y_t-z_t^{\splus}\rangle,
\end{align*}
where the equality follows from \eqref{DP:AL:F}  and \eqref{def:Psi:t},
the first inequality is due to Lemma~\ref{lem:eps.subgrad}, and
the last inequality is due to \eqref{eq:B_t.positive}. Using now the AM-GM inequality $a b \le (a^2 + b^2)/2$ with $(a,b)=(2\sqrt{2}\|w_t^{\splus}\|, (1/2\sqrt{2})\|z_t^{\splus} - y_t\|)$, and simple calculations we conclude from the previous inequality
%, using the Cauchy--Schwarz inequality and the , we obtain
that the following one holds:
\begin{align}
 {\cal L}_c(y_{<t},y_{t}, y_{>t}; q)&-{\cal L}_c(y_{<t}, z_t^{\splus}, y_{>t}; q) \nonumber\\ 
&\ge  \frac{5}{8\lam_t^{\splus}}\|z_t^{\splus}-y_t\|^2 + \frac{c}{4} \|A_t (z_t^{\splus} - y_t)\|^2  -\frac{1}{\lam_t^{\splus}}\left(\| 2\sqrt{2} w_t^{\splus}\| \left\|\frac1{2\sqrt{2}} (z_t^{\splus}-y_t)\right\|+ 2\epsilon_t^{\splus}\right) \nonumber\\
&\ge \frac{5}{8\lam_t^{\splus}}\|z_t^{\splus} - y_t\|^2 + \frac{c}{4} \|A_t (z_t^{\splus} - y_t)\|^2  -\frac{1}{\lam_t^{\splus}} \left(  4\|w_t^{\splus}\|^2 + \frac1{16} \|z_t^{\splus} - y_t\|^2+ 2\epsilon_t^{\splus} \right). \nonumber
\end{align}
The conclusion then follows by substituting the error bound in \eqref{eq:ref:lemma:A1} into the preceding inequality. 
\end{proof}

\begin{lemma}[Block-Sufficient Decrease]\label{lem:block.suff.decrease}
    Let $t \in \{1,\ldots,B\}$ be given and $\tau\in[0,1]$.
    If the prox stepsize $\lam_t^{\splus}$ satisfies $1/\lam_t^{\splus} \in [2m_t,\infty)$, then   \begin{align}\label{eq:lemma.block.decrease}
{\cal L}_c(y; q)-{\cal L}_c(y_{<t}, y_t(\tau), y_{>t}; q) & \geq \frac{\tau^2}{2\lambda_t^{\splus}}\|z_t^{\splus}-y_t\|^2 + \frac{c\tau}{4}\|A_t(z_t^{\splus} - y_t)\|^2,
\end{align}
where $y(\tau)$ is as in \eqref{def:z_tau}.
\end{lemma}

\begin{proof}
    We first claim that
\begin{align}\label{eq:diff_ztau2}
{\cal L}_c(y; q)-{\cal L}_c(y_{<t}, y_t(\tau), y_{>t}; q) &\geq \tau \Big[{\cal L}_c(y; q)-{\cal L}_c(y_{<t}, z_t^{\splus}, y_{>t}; q) \Big] + \frac{\tau^2-\tau}{2\lambda_t^{\splus}}\|z_t^{\splus}-y_t\|^2.
\end{align}
Indeed, since both $\Psi_t(\cdot)$ and $\lambda_t^{\splus}h_t(\cdot )$ are convex (due to Lemma \ref{lem:Psi.convexity}(b) and assumption (A1), respectively) and $\lambda_t^{\splus}>0$, the function $\Psi_t(\cdot)+\lambda_t^{\splus}h_t(\cdot)$ is convex. Moreover, using the identity $y_t(\tau)-y_t=\tau(z_t^{\splus} - y_t)$, we obtain
\[
(1-\tau)(\Psi_t+\lambda_t^{\splus}h_t)(y_t) + \tau(\Psi_t+\lambda_t^{\splus}h_t)(z_t^{\splus}) \geq (\Psi_t+\lambda_t^{\splus}h_t)(y_t(\tau)).
\]
Rearranging the above inequality and using \eqref{def:Psi:t} and the definition of ${\cal L}_c$ we obtain \eqref{eq:diff_ztau2}.

Finally, combining \eqref{eq:diff_ztau2} with \eqref{eq:lemma.block.decrease.hatz}, we conclude that \eqref{eq:lemma.block.decrease} holds.
\end{proof}

\begin{lemma}[Sufficient Decrease]\label{lem:suff.decrease}
For every $t \in \{1,\ldots,B\}$,  the $t$-th {\PABIPP} loop stops  in at most $1+\lceil \log(1+2m_t\lam_t) \rceil$ iterations and the pair $(y(\bar \tau),\lam^{\splus})$ satisfies
\begin{equation}\label{eq:renato.paper3}
 \Lc(y;q) - \Lc(y(\bar \tau);q) \geq \frac{\bar \tau^2}{2} \sum_{t=1}^B\frac{\|z_t^{\splus} - y_t\|^2}{\lambda_t^{\splus}} + \frac{c\bar \tau}{8}\sum_{t=1}^{B}\|A_t(z_t^{\splus}- y_t)\|^2,
\end{equation}
where $\bar \tau$ is as in line \ref{ln:PABIPP:tau} of {\PABIPP}.
\end{lemma}

\begin{proof}
Consider an arbitrary block index $t\in\{1,\ldots,B\}$. 
We first prove that the $t$-th {\PABIPP} loop stops  in at most $1+\lceil \log(1+2m_t\lam_t) \rceil$ iterations.
Firstly, we note that if the prox stepsize $\lam_t^{\splus}$ at a certain iteration of the $t$-th {\PABIPP} loop satisfies the condition $1/\lam_t^{\splus} \in [2m_t,\infty)$, then the Block-Sufficient Decrease Lemma (Lemma \ref{lem:block.suff.decrease}) guarantees that \eqref{eq:lemma.block.decrease} holds with $\tau=\bar \tau$. Hence, the inequality in \eqref{algo:PABIPP.adaptive.lambda} does not hold, i.e., the test in line~\ref{ln:PABIPP:test} of {\PABIPP} fails, and the $t$-th {\PABIPP} loop terminates at that iteration. On the other hand, let $\lambda_t>0$ be arbitrary. 
From the previous argument and the fact that $1/\lam_t^{\splus}$ is doubled each time line~\ref{ln:PABIPP:halve} of {\PABIPP} is executed, it follows that the $t$-th {\PABIPP} loop must terminate in at most $1+\lceil \log(1+2m_t\lam_t) \rceil$ iterations, for any choice of the initial proximal stepsize $\lam_t$.

We now prove that $(y(\bar \tau),\lambda^{\splus})$ satisfies \eqref{eq:renato.paper3}.
By the first part, for every $t\in\{1,\ldots,B\}$, the $t$-th {\PABIPP} loop terminates with a prox stepsize $\lambda_t^{\splus}$ for which the test in line~\ref{ln:PABIPP:test} of {\PABIPP} fails, i.e., the inequality in \eqref{algo:PABIPP.adaptive.lambda} does not hold. Since $y_t(\bar\tau)=y_t+\bar\tau(z_t^{\splus}-y_t)$ by \eqref{def:z_tau}, this means that
\begin{equation*}
\Lc(y;q) - \Lc(y_{<t},y_t(\bar\tau),y_{>t};q) \geq \frac{\bar\tau^2}{2\lambda_t^{\splus}}\|z_t^{\splus}-y_t\|^2 + \frac{c\bar\tau}{4}\|A_t(z_t^{\splus}-y_t)\|^2 \quad \forall t\in\{1,\ldots,B\}.
\end{equation*}
Note that this inequality is guaranteed by the test in line~\ref{ln:PABIPP:test} of {\PABIPP}, and hence holds even if the accepted stepsize does not satisfy $1/\lambda_t^{\splus}\geq 2m_t$. Combining the above inequality with Lemma~\ref{lem:jacobi.decrease} with $\tau=\bar\tau$, we obtain
\begin{align*}
\Lc(y;q) - \Lc(y(\bar\tau);q) &\geq \sum_{t=1}^{B}\Big[\Lc(y;q) - \Lc(y_{<t},y_{t}(\bar\tau),y_{>t};q)\Big]
- \frac{c(B-1)\bar\tau^2}{2}\sum_{t=1}^{B}\|A_t(z_t^{\splus}-y_t)\|^2\\
&\ge \sum_{t=1}^{B}\Big[\frac{\bar\tau^2}{2\lambda_t^{\splus}}\|z_t^{\splus}-y_t\|^2 + \frac{c\bar\tau}{4}\|A_t(z_t^{\splus} - y_t)\|^2\Big]- \frac{c(B-1)\bar\tau^2}{2}\sum_{t=1}^{B}\|A_t(z_t^{\splus}-y_t)\|^2\\
&=\frac{\bar\tau^2}{2}\sum_{t=1}^{B}\frac{\|z_t^{\splus}-y_t\|^2}{\lambda_t^{\splus}}+\frac{c\bar\tau[1-2(B-1)\bar\tau]}{4}\sum_{t=1}^{B}\|A_t(z_t^{\splus}-y_t)\|^2.
\end{align*}
Finally, substituting $\bar \tau=1/[4(B-1)]$, 
    simplifies the term $1 - 2(B-1)\bar\tau$ to $1/2$, directly yielding \eqref{eq:renato.paper3}.
\end{proof}

We now prove Proposition~\ref{prop:PABIPP}.

\noindent
{\bf Proof of Proposition \ref{prop:PABIPP}:}
(a) Let $t \in \{1,\ldots, B\}$ be given. Using the same argument as in Lemma~\ref{lem:suff.decrease}, we see that 
the $t$-th {\PABIPP} loop stops in at most $1+\lceil \log(1+2m_t\lam_t) \rceil$ iterations, the inequality in \eqref{eq:ABIPP:output.lambda.bound.2} holds, and \begin{equation*}
      \frac1{\lam_t^{\splus} } \le \max \left\{ \frac1{\lam_t} , 4m_t\right\}.
\end{equation*}
The identity in \eqref{eq:ABIPP:output.lambda.bound.2} now follows immediately from 
the inequality $1/\lam_t \le 1/\lam_t^{\splus}$ and the one above.

% Recall from Lemma~\ref{lemma:block-suff.decrease} that if  $1/\lambda_t^{\splus}\in[2m_t,\infty)$ 
% at some iteration of the $t$-th {\PABIPP} loop, then the loop must terminate at this iteration.
% Using this observation and the fact that $1/\lam_t^{\splus}$ is doubled every time line \ref{algo:ABIPP.if.loop} of {\PABIPP} is executed, we easily see that the $t$-th {\PABIPP} loop stops in at most $1+\lceil \log(1+2m_t\lam_t) \rceil$ iterations, the inequality in \eqref{eq:ABIPP:output.lambda.bound.2} holds, and \begin{equation*}
%       \frac1{\lam_t^{\splus} } \le \max \left\{ \frac1{\lam_t} , 4m_t\right\}.
% \end{equation*}
% The identity in \eqref{eq:ABIPP:output.lambda.bound.2} now follows immediately from 
% the inequality $1/\lam_t \le 1/\lam_t^{\splus}$ and the one above.

(b) First, we prove the inclusion in \eqref{prop.B-IPP.inclusion}.  
Note that, for each
$t\in\{1,\ldots,B\}$, the definition of ${\cal L}_{c}^{(\textrm{s})}(\cdot\,;q)$ in \eqref{def:smooth:ALM}, relation
\eqref{eq:ref:lemma:A1}, together with relation \eqref{prop1:subd} applied with $(\varepsilon, \beta)=(\epsilon_t^{\splus}, \lambda_t^{\splus})$,
gives the inclusion
\begin{align}\label{eq:COA1}
\frac{w_t^{\splus}}{\lambda_t^{\splus}}
&\overset{\eqref{eq:ref:lemma:A1}}\in \nabla_t  f(y_{<t},z_t^{\splus},y_{>t})+A_{t}^{*}\left[ q+c[A(y_{< t},z_t^{\splus},y_{>t})-b]\right]+\frac{1}{\lambda_t^{\splus}}(z_t^{\splus}-y_t)+\partial_{(\epsilon_t^{\splus}/\lambda_t^{\splus})}h_t(z_t^{\splus}). 
\end{align}
Now, it is straightforward to observe that
\[
A(y_{<t},z_t^{\splus},y_{>t})-b
=Az^{\splus}-b-\sum_{\substack{s=1, s\neq t}}^B
A_s(z_s^{\splus}-y_s).
\]
Thus, substituting this identity into \eqref{eq:COA1} yields
\[
\frac{w_t^{\splus}}{\lambda_t^{\splus}}\in \nabla_t  f(y_{<t},z_t^{\splus},y_{>t})+A_{t}^{*}\left(q+c(Az^{\splus}-b) - c\sum_{s=1,s\neq t}^{B}A_s (z_s^{\splus}-y_s)\right)+\frac{1}{\lambda_t^{\splus}}(z_t^{\splus} -y_t)+\partial_{(\epsilon_t^{\splus}/\lambda_t^{\splus})}h_t(z_t^{\splus}).
\]

%{\bf Working ...}

Hence, using the separability of $f$, which implies that $\nabla_t  f(y_{<t},z_t^{\splus},y_{>t})=\nabla f_t(z_t^{\splus})$,
and recalling the definition of $r_t^{\splus}$ in line \ref{ln:PABIPP:residual} of {\PABIPP}, we conclude that, for every $t\in \{1,\ldots,B\}$, the following inclusion holds:
% that for every $t\in \{1,\ldots,B\}$,
\[
r_t^{\splus} \in \nabla f_t(z_t^{\splus})+\partial_{(\epsilon_t^{\splus}/\lambda_t^{\splus})}h_t(z_t^{\splus})+A_t^*\left[q+c(Az^{\splus}-b)\right]. %\quad \forall t\in \{1,\ldots,B\}.
\]
On the other hand, \eqref{prop2:subd} and the definition of
$\delta^{\splus}$ in line \ref{ln:PABIPP:assemble}  
of {\PABIPP}
imply
that 
\[
\partial_{(\epsilon_1^{\splus}/\lambda_1^{\splus})}h_1(z_1^{\splus}) \times \ldots \times
\partial_{(\epsilon_B^{\splus}/\lambda_B^{\splus})}h_B(z_B^{\splus})  \subseteq \partial_{\delta^{\splus}}h(z^{\splus}).
\]
Thus, we conclude that \eqref{prop.B-IPP.inclusion} holds by combining the last two inclusions and using the definition of $r^{\splus}$ in line \ref{ln:PABIPP:assemble} of {\PABIPP}.

It remains to verify \eqref{prop.B-IPP.inequality}.
First, we 
set
$
\Delta{\cal L}_c
:=
{\cal L}_c(y;q)-{\cal L}_c(y^{\splus};q),$
where
$
y^{\splus}=y(\bar\tau).
$
Then, the Sufficient Decrease Lemma (Lemma~\ref{lem:suff.decrease}) implies that
\begin{equation}\label{eq:adapt.Delta.Lc.ineq}
\Delta {\cal L}_c  \ge \frac{\bar \tau^2}{2}\sum_{t =1}^B \frac{\|z_t^{\splus}-y_t\|^2}{\lam_t^{\splus}} + \frac{c\bar \tau}{8}\sum_{t =1}^B \|A_t (z_t^{\splus}-y_t)\|^2.
\end{equation}
%where $y^{\splus}=y(\bar \tau)$.

Next, we estimate the residuals   
$r^{\splus}$ and $\delta^{\splus}$, which depend on the block residuals $w_t^{\splus}$ and $\epsilon_t^{\splus}$.  
Thus, we estimate the former block residuals by 
using \eqref{eq:ref:lemma:A1}, \eqref{eq:adapt.Delta.Lc.ineq},  and that $1/\lam_t^{\splus}\leq 1/\min(\lambda^{\splus})$ (see \eqref{def.lambda.min.max}), that is,
\begin{align*}
\sum_{t=1}^B \left( 2\frac{\|w_t^{\splus}\|^2}{(\lam_t^{\splus})^2}  + \frac{\epsilon_t^{\splus}}{\lam_t^{\splus}} \right)
&\le \left[\frac{2}{\min(\lambda^{\splus})}+1\right] \sum_{t=1}^B \left( \frac{\|w_t^{\splus}\|^2 + \epsilon_t^{\splus}}{\lam_t^{\splus}}\right)\\
&\overset{\eqref{eq:ref:lemma:A1}}\le \left[\frac{2}{\min(\lambda^{\splus})}+1\right]  \sum_{t=1}^B \left( \frac{\|z_t^{\splus}-y_t\|^2}{32\lam_t^{\splus}}  \right) \overset{\eqref{eq:adapt.Delta.Lc.ineq}}\le  \frac{1}{16\bar\tau^2}\left[\frac{2}{\min(\lambda^{\splus})}+1\right] \Delta {\cal L}_c.
\end{align*}
Now, setting $D_{t}:=\|r_t^{\splus}-w_t^{\splus}/\lam_t^{\splus}\|^2$, using that $\|a+b\|^2\leq 2\|a\|^2+2\|b\|^2$, for any $a,b\in \R^n$,
together with the definition of $\delta^{\splus}$ (see line \ref{ln:PABIPP:assemble} of {\PABIPP}) gives
\begin{align}\label{eq:suff.decr.inter}
\|r^{\splus}\|^2 &+\delta^{\splus} =
\sum_{t=1}^B \left( \|r_t^{\splus}\|^2 + \frac{\epsilon_t^{\splus}}{\lam_t^{\splus}}  \right) 
\le 
\sum_{t=1}^B \left( 2 D_{t} + 2\frac{\|w_t^{\splus}\|^2}{(\lam_t^{\splus})^2}  + \frac{\epsilon_t^{\splus}}{\lam_t^{\splus}}
\right) \leq 2 \sum_{t=1}^B  D_{t} +   \frac{1}{16\bar\tau^2}\left[\frac{2}{\min(\lambda^{\splus})}+1\right] \Delta {\cal L}_c.
\end{align}
We now control the first term on the right-hand side of the previous inequality.

Using the definition of
$r_t^{\splus}$ in line \ref{ln:PABIPP:residual} of {\PABIPP}, we have
\[
r_t^{\splus}-\frac{w_t^{\splus}}{\lambda_t^{\splus}}
=
cA_t^*
\left[
\sum_{\substack{s=1, s\neq t}}^B
A_s(z_s^{\splus}-y_s)
\right]
-
\frac{z_t^{\splus}-y_t}{\lambda_t^{\splus}}.
\]
Consequently, it follows from the definition of $D_t$ and simple calculations that
\begin{align*}
D_t
&\leq
2c^2\|A_t\|^2
\left(
\sum_{\substack{s=1, s\neq t}}^B
\|A_s(z_s^{\splus}-y_s)\|
\right)^2
+
2\frac{\|z_t^{\splus}-y_t\|^2}
{(\lambda_t^{\splus})^2}
\\
&\leq
2c^2\|A_t\|^2(B-1)
\sum_{\substack{s=1, s\neq t}}^B
\|A_s(z_s^{\splus}-y_s)\|^2
+
\frac{2}{\min(\lambda^{\splus})}
\frac{\|z_t^{\splus}-y_t\|^2}
{\lambda_t^{\splus}},
\end{align*}
where the second inequality follows from Cauchy--Schwarz.
Applying
\eqref{eq:adapt.Delta.Lc.ineq} to the first term, and using the fact that
$B-1=1/(4\bar\tau)$ by the definition of $\bar\tau$ in line~\ref{ln:PABIPP:tau} of {\PABIPP}, yields
\begin{align}
D_t
\leq
2\left[
2c\|A_t\|^2\bar\tau^{-2}\Delta{\cal L}_c
+
\frac{1}{\min(\lambda^{\splus})}
\frac{\|z_t^{\splus}-y_t\|^2}
{\lambda_t^{\splus}}
\right].
\label{eq:bound:vi:first:Ine}
\end{align}

Summing \eqref{eq:bound:vi:first:Ine} over
$t=1,\ldots,B$, using the definition of
$\|A\|_\dagger$ in \eqref{eq:def.zeta2andblock_norm}, and invoking
\eqref{eq:adapt.Delta.Lc.ineq} once more, we obtain
\begin{align*}
\sum_{t=1}^B D_t
&\leq
4c\|A\|_\dagger^2\bar\tau^{-2}\Delta{\cal L}_c
+
\frac{2}{\min(\lambda^{\splus})}
\sum_{t=1}^B
\frac{\|z_t^{\splus}-y_t\|^2}
{\lambda_t^{\splus}}
\leq
4\bar\tau^{-2}
\left[
c\|A\|_\dagger^2
+
\frac{1}{\min(\lambda^{\splus})}
\right]
\Delta{\cal L}_c.
\end{align*}
Now, using this bound,  \eqref{eq:suff.decr.inter}, the
definition of $\zeta$ in \eqref{eq:def.zeta2andblock_norm}, gives precisely
\eqref{prop.B-IPP.inequality}. This completes the proof. \qedflush

\section{Proof of Theorem~\ref{thm:main} }
\label{sec:proof.thm}

This section proves Theorem~\ref{thm:main}. The core of the proof is the analysis of a single stage of {\PAADMM}, which is done in Proposition~\ref{prop:stage.complexity} at the end of this section. In summary, this proposition shows that each stage ends with a $\rho$-stationary solution of \eqref{initial.problem}--\eqref{eq:block_structure} after a number of inner iterations bounded in terms of the stage input $(x^{s-1},p^{s-1},c_{s-1})$. It also shows that the Lagrange multiplier and the scaled feasibility residual satisfy local bounds.

The proof of Proposition~\ref{prop:stage.complexity} relies on Lemmas~\ref{lem:residual.potential}--\ref{lem:mult.variation}. In particular, Lemmas~\ref{lem:residual.potential}--\ref{lem:AL.gap} relate the stationarity residual to the potential $T_i$ (see line~\ref{ln:PAADMM:T} of {\PAADMM}) and bound $T_i$. Lemmas~\ref{lem:mult.step} and~\ref{lem:mult.bounds} bound the Lagrange multipliers and the scaled feasibility residual, and Lemma~\ref{lem:mult.variation} uses the LMT to bound the cumulative variation of the multipliers.

Before proving Theorem~\ref{thm:main}, we introduce the necessary notation and provide basic facts about {\PAADMM} that will be used throughout the subsequent analysis.

In this part of the analysis, we fix an arbitrary outer-loop index $s$ generated by {\PAADMM}. Hence, the penalty parameter $c=c_{s-1}$ is kept fixed, the stage starts from $(y^0,q^0,\lambda^0) = (x^{s-1},p^{s-1},\gamma^{s-1})$, and the sequences
\[
\{(y^i,q^i,\lambda^i,T_i,R_i)\}_{i\geq 0}, \qquad \{(z^{i}, r^{i},\delta_{i})\}_{i\geq 1},
\]
refer to the (inner) quantities generated during this fixed \textit{stage}. Note that the second sequence is indexed from $i=1$, since its elements are only computed from the first inner iteration onward, whereas the first one also includes the initial quantities. For every $k \ge 1$, define the $k$-th epoch $\calI_k$
as the index set 
\begin{equation*}
{\calI}_k := \{\cali_{k-1}+1, \ldots, \cali_k\},
\end{equation*}
where $\cali_k$ is computed either in line~\ref{ln:PAADMM:count.stop} or line~\ref{ln:PAADMM:count.LMT} of {\PAADMM} (recall from the second paragraph following {\PAADMM} that $\cali_k$ denotes the inner iteration at which the $k$-th Lagrange multiplier update occurs), with the convention that $\cali_0=0$.
Moreover, let
\begin{equation}\label{def.variablesTilde}
   (\ty^k,\tq^k,\tT_k):=(y^{\cali_k},q^{\cali_k},T_{\cali_k}) \quad \forall k\geq 0.
\end{equation}

We now make three additional remarks  about the logic of {\PAADMM} regarding the prox stepsize and the Lagrange multiplier. First, due to the definition of $\cali_k$, it follows that $q^i=q^{i-1}$ for every $i \in \{\cali_{k-1}+1,\ldots, \cali_k-1\}=\calI_k\setminus \{\cali_k\}$, which implies that 
\begin{equation}\label{eq:q.constant}
    q^{i-1}=q^{\cali_{k-1}}=\tq^{k-1} \ \ \forall i\in \calI_k.
\end{equation}
Moreover, \eqref{def.variablesTilde} and \eqref{eq:q.constant} with $i=\cali_k$
imply that 
\begin{equation}\label{eq:q.update}
    \tq^{k} = q^{\cali_k}= q^{\cali_k-1} + c(Az^{\cali_k}-b)=\tq^{k-1} +  c(Az^{\cali_k}-b) \quad \forall k\geq 1.
\end{equation}
Second, $\tilde q^0 \in A(\R^n)$ due to  the facts that  $\tilde q^0=q^0$
by \eqref{def.variablesTilde} with $k=0$,
and the input $q^0$ of  {\PAADMM} is in $A(\R^n)$. 
This observation, 
the fact that $b \in A(\R^n)$ by (A3), identity \eqref{eq:q.update}, and
a simple induction argument, then imply that
\begin{equation}\label{eq:q.in.Img.A}
   \tilde q^k \in A(\R^n) \quad \forall  k \ge 0, 
\end{equation}
and hence that
\begin{equation}\label{eq:q.in.Img.A.2}
   q^i \in A(\R^n) \quad \forall  i \ge 0,
\end{equation}
in view of~\eqref{eq:q.constant}.

Third, we clarify that the quantity \(R_i\) is well defined for every \(i\geq0\). For \(i=0\), it is set to one in line~\ref{ln:PAADMM:stage.init} of {\PAADMM}. For \(i\geq1\), it is computed in line~\ref{ln:PAADMM:R} of {\PAADMM}, and we adopt the conventions that the maximum over an empty set is zero and that the ratio defining \(R_i\) is equal to one whenever both its numerator and denominator vanish. Indeed, the denominator vanishes if and only if
\[
q^j=q^{j-1}\quad \forall j=1,\ldots,i-1,
\qquad\text{and}\qquad
\mathsf q^{\mathrm{cand}}=q^{i-1},
\]
in which case the numerator also vanishes. Moreover, these conventions ensure that $1\leq R_i\leq \max\{1,i\}$ for all $i\geq 0$.

With the necessary notation now in place, we establish the intermediate results needed to prove Proposition~\ref{prop:stage.complexity}.

The next lemma shows that for every inner iteration $i$, the quadruple $(z^i,q^{i-1}+c(Az^i-b),r^i,\delta_i)$ satisfies the inclusion in~\eqref{eq:rho.eta.solution}, where the quintuple $(y^i,\lambda^i;z^i,r^i,\delta_i)$ is generated from {\PABIPP} (see line~\ref{ln:PAADMM:call} of {\PAADMM}). The lemma also provides a key bound relating the stationarity residual $\|r^i\|^2+\delta_i$ to the potential increment $T_i-T_{i-1}$.

\begin{lemma}\label{lem:residual.potential}
    Let $s$ be a stage generated by {\PAADMM} and let $i$ be an iteration within this stage. Then, it holds that
    \begin{align}
r^i \in \nabla f(z^i)+&\partial_{\delta_i}h(z^i)
+A^*\left[q^{i-1}+c(Az^i-b)\right],
\label{eq:C1a}
\\
\|r^i\|^2+\delta_i
&\le
\left(\tilde\zeta + c \zeta\right)(T_i-T_{i-1}),
\label{eq:C1b}
\end{align}
where $\zeta$ and $\tilde \zeta$ are as in \eqref{eq:def.zeta2andblock_norm} and \eqref{eq:tildezeta.Upsilon}, respectively.
\end{lemma}

\begin{proof}
It follows from lines \ref{ln:PAADMM:call} and \ref{ln:PAADMM:T} of {\PAADMM}, and Proposition~\ref{prop:PABIPP} with $(y,q,\lam,c)=(y^{i-1},q^{i-1},\lam^{i-1},c)$
and $(y^{\splus},\lam^{\splus};z^{\splus},r^{\splus},\delta^{\splus}) = (y^i,\lam^i;z^i,r^i,\delta_i)$,
that the inclusion \eqref{eq:C1a} holds and
\begin{align}\label{lemma:norm:residual:Ine.2:sec:6.2:2}
\|r^i\|^2 +\delta_i  
\le \left[\frac{1}{16\bar\tau^2}\left(\frac{130}{\min(\lambda^{i})}+1\right) + c \zeta\right](T_i-T_{i-1}),
\end{align}
\begin{equation}\label{cond:lam:sec:6.2}
\lam_t^i \le \lam_t^{i-1}, \quad \max\left\{ \frac{1}{\lam_t^{i}},4m_t\right\} = \max\left\{ \frac{1}{\lam_t^{i-1}},4m_t\right\} \quad \forall t \in \{1,\ldots,B\},
\end{equation}
% \begin{equation}\label{lemma:inclusion:residual:Ine.2:sec:6.2}
% r^i  \in \nabla f (u^i)+ \partial_{\delta_i} h(u^i) +
% A^*[q^{i-1}+c(Au^i-b)],
% \end{equation}
where $\zeta$ is as in \eqref{eq:def.zeta2andblock_norm}.
%Hence, \eqref{eq:C1a} follows immediately from~\eqref{lemma:inclusion:residual:Ine.2:sec:6.2}.
We now use \eqref{lemma:norm:residual:Ine.2:sec:6.2:2} and \eqref{cond:lam:sec:6.2} to prove \eqref{eq:C1b}. Indeed, it follows from \eqref{cond:lam:sec:6.2} and the way $\gamma^0,\ldots,\gamma^{s-1}$ are chosen in line \ref{ln:PAADMM:stage.output} of {\PAADMM} that
\begin{equation*}
\frac{1}{\lam_t^i} \leq \max\left\{ \frac{1}{\lam_t^{i}},4m_t\right\} = \max\left\{ \frac{1}{\gamma^{0}_t},4m_t\right\} \leq \frac{1}{\gamma_t^0} + 4m_t \quad \forall t \in \{1,\ldots,B\}.
\end{equation*}
The last conclusion, the definition of $\tilde \zeta$ as in \eqref{eq:tildezeta.Upsilon}, and  inequality \eqref{lemma:norm:residual:Ine.2:sec:6.2:2} now immediately imply \eqref{eq:C1b}.
\end{proof}

Summing \eqref{eq:C1b} from $j=1$ to $j=i$ and using the fact that $T_0=0$, we have that
\[
\min_{1\leq j\leq i}(\|r^j\|^2+\delta_j) \leq \frac{1}{i}\sum_{j=1}^{i}(\|r^j\|^2+\delta_j) \leq \frac{(\tilde \zeta + c\zeta)T_i}{i}.
\]
Consequently, it follows from the (inner) termination criteria in line \ref{ln:PAADMM:stat.test} of {\PAADMM} that an upper bound on $T_i$ will immediately yield an upper bound on the total number of inner iterations. Hence, the next two lemmas provide an intermediate step toward this goal. Specifically, Lemma~\ref{lem:Ti.identity} provides the exact expression for $T_i$, and Lemma~\ref{lem:Ti.bound} establishes a bound to that expression.

\begin{lemma}\label{lem:Ti.identity}
Let $k$ be an epoch generated by the $s$-th stage of
{\PAADMM}. Then, for every $i\in\calI_k$,
\begin{equation}\label{eq:Ti.value}
T_i = \left [ {\cal L}_c(y^{0};q^{0}) - {\cal L}_c(y^i;\tq^{k-1}) \right ]+(1-\bar\tau)\sum_{j=1}^{i-1} \Inner{A(y^{j-1}- z^{j})}{q^{j}-q^{j-1}}+\frac{1}{c}\sum_{j=1}^{i-1}\|q^j-q^{j-1}\|^2,
\end{equation} 
where $\bar \tau$ is as in line \ref{ln:PABIPP:tau} of {\PABIPP}.
\end{lemma}

\begin{proof}
Assume that $i$ is an inner-iteration index generated by {\PAADMM} and that $i\in{\calI}_k$.
By the definition of $T_i$ as in line \ref{ln:PAADMM:T} of {\PAADMM}, we have
\begin{align}
\label{eq.Ti.sum.1}
T_i-T_1
&=\sum_{j=1}^{i-1}[T_{j+1}-T_{j}]=
\sum_{j=1}^{i-1}
\left[
{\cal L}_c(y^j;q^j)
-
{\cal L}_c(y^{j+1};q^j)
\right]   \nonumber \\
&=
\sum_{j=1}^{i-1}
\left[
{\cal L}_c(y^j;q^j)
-
{\cal L}_c(y^j;q^{j-1})
\right]                                   +
\sum_{j=1}^{i-1}
\left[
{\cal L}_c(y^j;q^{j-1})
-
{\cal L}_c(y^{j+1};q^j)
\right].
\end{align}
Since $i\in{\calI}_k$, we have
$q^{i-1}=\tq^{k-1}$ by~\eqref{eq:q.constant}, and hence,
\begin{align*}
T_1
+
\sum_{j=1}^{i-1}
\left[
{\cal L}_c(y^j;q^{j-1})
-
{\cal L}_c(y^{j+1};q^j)
\right]
=
T_1+
{\cal L}_c(y^1;q^0)
-
{\cal L}_c(y^i;q^{i-1})                   &=
{\cal L}_c(y^0;q^0)
-
{\cal L}_c(y^i;\tq^{k-1}).
\end{align*}
Combining~\eqref{eq.Ti.sum.1} with the previous identity, we obtain
\begin{equation}
\label{eq:Ti.new.expression}
T_i
=
\sum_{j=1}^{i-1}
\left[
{\cal L}_c(y^j;q^j)
-
{\cal L}_c(y^j;q^{j-1})
\right]
+
{\cal L}_c(y^0;q^0)
-
{\cal L}_c(y^i;\tq^{k-1}).
\end{equation}

We now calculate the first term on the right-hand side of~\eqref{eq:Ti.new.expression}. We claim that, for every
$j\in\{1,\ldots,i-1\}$,
\begin{align}
\label{eq:proof.Ti.part2}
\Lc(y^{j};q^{j})-\Lc(y^{j};q^{j-1})
=
(1-\bar \tau)\Inner{A(y^{j-1}-z^j)}{q^j-q^{j-1}}
+
\frac{1}{c}\|q^j-q^{j-1}\|^2 .
\end{align}
Indeed, if $j\in{\calI}_\ell\setminus\{\cali_\ell\}$ for some epoch index $\ell\in\{1,\ldots,k\}$, then $q^j=q^{j-1}$ due to the remark containing \eqref{eq:q.constant}, and hence both sides of~\eqref{eq:proof.Ti.part2} are equal to zero. 
It remains to consider the case $j=\cali_\ell$ for some epoch index $\ell$. Using the definition of the augmented
Lagrangian as in \eqref{DP:AL:F}, we have
\begin{align*}
{\cal L}_c(y^j;q^j)-{\cal L}_c(y^j;q^{j-1})&={\cal L}_c(y^{\cali_\ell};q^{\cali_\ell})-{\cal L}_c(y^{\cali_\ell};q^{\cali_\ell-1}) \\
&=
\Inner{Ay^{\cali_\ell}-b}{q^{\cali_\ell}-q^{\cali_\ell-1}}\\
&=\Inner{A(y^{\cali_\ell}-z^{\cali_\ell})+(Az^{\cali_\ell}-b)}{q^{\cali_\ell}-q^{\cali_\ell-1}
} \\
&=
\Inner{A(y^{\cali_\ell}-z^{\cali_\ell})}{q^{\cali_\ell}-q^{\cali_\ell-1}}+\Inner{Az^{\cali_\ell}-b
}{q^{\cali_\ell}-q^{\cali_\ell-1}}.
\end{align*}
Since the multiplier update gives $q^{\cali_\ell}-q^{\cali_\ell-1}=c(Az^{\cali_\ell}-b)$, and $y^{\cali_\ell} = (1-\bar \tau)y^{\cali_\ell-1}+\bar\tau z^{\cali_\ell}$, we conclude that \eqref{eq:proof.Ti.part2} holds, and \eqref{eq:Ti.value} follows by combining~\eqref{eq:Ti.new.expression} and  \eqref{eq:proof.Ti.part2}.
\end{proof}

\begin{lemma}\label{lem:Ti.bound}
Let $k$ be an epoch generated by the $s$-th stage of
{\PAADMM}. Then, for every $i\in\calI_k$,
\begin{equation}
\label{eq:Ti.bound.final}
T_i
\le
2\left[
{\cal L}_c(\ty^0;\tq^0)
-
{\cal L}_c(y^i;\tq^{k-1})
\right]
+
\left[\frac{(1-\bar\tau)}{\chi}+\frac{2}{c}\right]
\sum_{\ell=1}^{k-1}
\|\tq^\ell-\tq^{\ell-1}\|^2 ,
\end{equation}
where $\bar \tau$ is as in line \ref{ln:PABIPP:tau} of {\PABIPP} and 
\begin{equation}\label{def:chi}
\chi:=\frac{c\bar\tau}{8B(1-\bar\tau)}.
\end{equation}
\end{lemma}

\begin{proof}
We use the expression for $T_i$ obtained in \eqref{eq:Ti.value}. We first bound the inner-product term appearing in that expression. Using the Cauchy--Schwarz inequality, followed by the AM--GM inequality, the fact that $\chi>0$, and the block norm inequality, we obtain, for $j=1,\ldots,i-1$,
\begin{align*}
 \Inner{A(y^{j-1}-z^j)}{q^j-q^{j-1}}
 &\leq
 \|A(y^{j-1}-z^j)\|\|q^j-q^{j-1}\|
 \\
 &\leq
 \frac{1}{2}\left[\chi \|A(y^{j-1}-z^j)\|^2 + \frac{1}{\chi}\|q^j-q^{j-1}\|^2\right]\\
 &\leq \frac{B\chi}{2}\sum_{t=1}^B\|A_t(y_t^{j-1}-z^j_t)\|^2 + \frac{1}{2\chi}\|q^j-q^{j-1}\|^2 
\end{align*}
where the last inequality follows from
$
\|w_1+\ldots+w_B\|^2\le B \left( \|w_1\|^2+\ldots+\|w_B\|^2 \right)
$
with $w_t = A_t(y_t^{j-1}-z_t^j)$ for $t=1,\ldots,B$. Now, using the definition of $\chi$ in \eqref{def:chi}, the previous inequality, Lemma~\ref{lem:suff.decrease} with $(z_t^{\splus},y_t,y,q,y(\bar\tau))=(z_t^j,y_t^{j-1},y^{j-1},q^{j-1},y^j)$, and  line \ref{ln:PAADMM:T} of {\PAADMM}, we conclude that
\begin{align*}
    (1-\bar\tau)\Inner{A(y^{j-1}-z^j)}{q^j-q^{j-1}}
 &\leq \frac{c\bar\tau }{16}\sum_{t=1}^B\|A_t(y_t^{j-1}-z^j_t)\|^2 + \frac{(1-\bar\tau)}{2\chi}\|q^j-q^{j-1}\|^2\\
 &\leq \frac{1}{2} \left[
{\cal L}_c(y^{j-1};q^{j-1})-{\cal L}_c(y^j;q^{j-1})
\right]+\frac{(1-\bar\tau)}{2\chi}
\|q^j-q^{j-1}\|^2 \\
& = \frac{1}{2}[T_j-T_{j-1}] + \frac{(1-\bar\tau)}{2\chi}
\|q^j-q^{j-1}\|^2.
\end{align*}
Summing the above inequality from $j=1$ to $j=i-1$, using the expression for $T_i$ in \eqref{eq:Ti.value},
and the fact that $T_0=0$, we have
\begin{align*}
T_i
&\le\left [ {\cal L}_c(y^{0};q^{0}) - {\cal L}_c(y^i;\tq^{k-1}) \right ]+\frac12 T_{i-1}+\left[\frac{(1-\bar\tau)}{2\chi}+\frac{1}{c}\right]\sum_{j=1}^{i-1}\|q^j-q^{j-1}\|^2 \nonumber\\
&=\left [ {\cal L}_c(y^{0};q^{0}) - {\cal L}_c(y^i;\tq^{k-1}) \right ]+\frac12 T_{i-1}+\left[\frac{(1-\bar\tau)}{2\chi}+\frac{1}{c}\right]\sum_{\ell=1}^{k-1}\|\tq^\ell-\tq^{\ell-1}\|^2,
\end{align*}
where the above equality follows from the fact that relation \eqref{def.variablesTilde}  and the first remark in the paragraph following
\eqref{def.variablesTilde}, both with $k=\ell$,
imply that
\begin{equation}
\label{eq:q-update}
q^j-q^{j-1}
=
\begin{cases}
0,
& \text{if } j\in \calI_{\ell}\setminus\{\cali_\ell\},\\
\tq^\ell-\tq^{\ell-1},
& \text{if } j=\cali_\ell
\end{cases}
\end{equation}
for every
 $1\leq \ell\leq k$.
The conclusion of the lemma now follows from the above inequality and the fact that $T_{i-1}\leq T_i$ (see \eqref{eq:C1b}).
\end{proof}

Lemma~\ref{lem:Ti.bound} does not yet provide the desired uniform stagewise bound. Indeed, if $i\in\calI_k$, then the right-hand side of~\eqref{eq:Ti.bound.final} still depends on the generated point $y^i$ and the multipliers $\{\tq^\ell\}_{\ell=0}^{k-1}$. More specifically, this dependence appears through an augmented-Lagrangian gap and the cumulative squared variation of the multipliers. We control these two quantities separately. Lemma~\ref{lem:AL.gap} provides a general bound for the augmented-Lagrangian gap. Lemmas~\ref{lem:mult.step} and~\ref{lem:mult.bounds} establish uniform bounds on the multiplier sequence and its individual increments. Finally, Lemma~\ref{lem:mult.variation} controls their cumulative variation. Together, these results allow us to express the upper bound on $T_i$ solely in terms of the problem parameters and the current stage's initial data.

\begin{lemma}\label{lem:AL.gap}
For any given 
$c >0$ and pairs
$(x,p) \in  {\cal H} \times \mathbb{R}^l$ and
$(y,q) \in  {\cal H} \times \mathbb{R}^l$, we have
\begin{equation}\label{def:diff:lagrangian}
{\cal L}_{c}(x; p)
- {\cal L}_{c}(y; q )
\le \phi_{\mathrm{sup}}-\phi_{\mathrm{inf}}
+ c\|Ax-b\|^2 +\frac{\|p\|^2+\|q\|^2}{2c}
\end{equation}
where  $(\phi_{\mathrm{sup}},\phi_{\mathrm{inf}} )$ is as in \eqref{def:damH}.
\end{lemma}
\begin{proof}
Using the definitions of $\phi_{\inf}$ and ${\cal L}_c(\cdot\,;\,\cdot)$ in \eqref{def:damH} and \eqref{DP:AL:F}, respectively, we have
\begin{align*}
{\cal L}_{c}(y ; q) - \phi_{\inf} &\overset{\eqref{def:damH}}\geq  {\cal L}_{c}(y ; q)-(f+h)(y) \\
&\overset{\eqref{DP:AL:F}}=\langle q, A y-b\rangle +\frac{c}{2}\|A y-b\|^2=\frac{1}{2}\left\|\frac{q}{\sqrt{c}}+\sqrt{c}(A y-b)\right\|^2-\frac{\|q\|^2}{2c} \ge -\frac{\|q \|^2}{2c}.
\end{align*}
On the other hand, using the definition of $\phi_{\sup}$ in \eqref{def:damH} and the Cauchy--Schwarz and Young inequalities, we have
\begin{align*}
{\cal L}_{c}(x; p) - \phi_{\sup} 
&\overset{\eqref{def:damH}}\le {\cal L}_{c}(x; p)-(f+h)(x)  \overset{\eqref{DP:AL:F}}= \langle p, A x -b\rangle +\frac{c\|A x-b\|^2}{2}\\
&\le \left( \frac{\|p\|^2}{2c}+\frac{c\|A x-b\|^2}{2} \right) +\frac{c\|A x-b\|^2}{2}=\frac{\|p\|^2}{2c}+c\|A x-b\|^2.
\end{align*}
 Combining the above two relations, we then conclude that \eqref{def:diff:lagrangian} holds.
\end{proof}

\begin{lemma}\label{lem:mult.step}
Let $s$ be a stage generated by {\PAADMM} and let
     ${\calI}_k$ be an epoch within its $s$-th stage  and define
\begin{equation}\label{def:I:k:C1}
    {\calI}_k(C) := \{i\in {\calI}_k:\|r^{i}\|^2+\delta_i\leq C^2\},
\end{equation}
where $C>0$ is part of the input for {\PAADMM} and $(r^i, \delta_i)$ is as in line~\ref{ln:PAADMM:call} of {\PAADMM}.
Then, the pair $(\tq^{k-1}, z^i)$ satisfies
\begin{equation}\label{eq:FeasBound.proof.1}
\|\tilde q^{k-1} + c_{s-1} (Az^i-b)\| \le 
\max \{ \|\tilde q^{k-1}\|,\Upsilon(C)\}, \quad \forall i \in \calI_k(C),
\end{equation}
where $\Upsilon(C)$ is as in \eqref{eq:tildezeta.Upsilon}. As a consequence,
\begin{equation}\label{step:induc:next:result}
  \|\tq^k\|\le  \max \{ \|\tilde q^{k-1}\|,\Upsilon(C)\}.
\end{equation}
\end{lemma}
\begin{proof}
We first argue that \eqref{step:induc:next:result} follows as an immediate consequence of \eqref{eq:FeasBound.proof.1}.
Indeed, noting that the assumption that $C \ge \rho$ and the logic of {\PAADMM} imply that $\cali_k$ always satisfies $\|r^{\cali_k}\|^2+\delta_{\cali_k}\leq C^2$, i.e., $\cali_k \in \calI_k(C)$, and using \eqref{eq:FeasBound.proof.1} with $i=\cali_k$, and identities \eqref{def.variablesTilde} and \eqref{eq:q.update}, we immediately see that \eqref{step:induc:next:result} holds.

To show \eqref{eq:FeasBound.proof.1}, let $i \in \calI_k(C)$ be given. The proof of \eqref{eq:FeasBound.proof.1} relies on Lemma~\ref{lem:mult.bound.aux} with $(q^-,\varrho) = (\tilde q^{k-1},c)$ and $(z,q,r,\delta)=(z^i, \underline{q}^i, \underline{r}^i, \delta_i)$, where
\begin{align}\label{eq:IkC.bounded.q.r}
\underline{q}^i:=\tilde q^{k-1}+c(Az^i-b) \quad \text{and} \quad \underline{r}^i := r^i-\nabla f(z^i).
\end{align}
We first claim that  $(q^-,\varrho) = (\tilde q^{k-1},c)$ and $(z,q,r,\delta)=(z^i, \underline{q}^i, \underline{r}^i, \delta_i)$ satisfy the assumptions of Lemma~\ref{lem:mult.bound.aux}, i.e., the inclusion $\tilde q^{k-1} \in A(\R^n)$ and the conditions in \eqref{eq:cond:Lem:A1}. Indeed, the inclusion is due to \eqref{eq:q.in.Img.A}. Moreover, the identity in \eqref{eq:cond:Lem:A1} is due to the first identity in \eqref{eq:IkC.bounded.q.r}. Also,  the inclusion~\eqref{eq:C1a} of Lemma~\ref{lem:residual.potential}, and
relations \eqref{eq:q.constant} and \eqref{eq:IkC.bounded.q.r}, imply that
$
\underline{r}^i \in  \partial_{\delta_i} h(z^i) + A^*\underline{q}^i
$. Thus, the claim holds.

Now, 
noting that the triangle inequality for norms, the second identity in \eqref{eq:IkC.bounded.q.r}, the definitions of $\Ik(C)$ in \eqref{def:I:k:C1} and $\widehat \nabla_f$ in \eqref{def:damH}, imply that
\begin{equation}\label{eq:ri.deltai.C}
   \delta_i + \|\underline{r}^i\| \overset{\eqref{eq:IkC.bounded.q.r}}{\leq} \delta_i + \|r^i\|+\|\nabla f(z^i)\| \overset{\eqref{def:I:k:C1}}{\le} C^2 + (C + \|\nabla f(z^i)\|) \overset{\eqref{def:damH}}{\le} C^2 + (C + \widehat \nabla_f),
\end{equation}
and using the conclusion of Lemma~\ref{lem:mult.bound.aux}, we conclude that
\begin{align*}
 \|\underline{q}^i\| &\overset{\eqref{q-bound-2}}{\leq} \max\left\{\| \tilde q^{k-1}\|\, , \,\Xi(\|z^i-\mathrm{\bar x}\|\, ,\, \|\underline{r}^i\|+\delta_i) \right\}\\
 &\overset{\eqref{eq:ri.deltai.C}}{\leq} \max\left\{\| \tilde q^{k-1}\|\, , \, \Xi(D_h \, , \, C+C^2 + \widehat \nabla_f) \right\} \leq \max\{\|\tq^{k-1}\|\, , \,\Upsilon(C)\},
\end{align*}
where the second inequality is due to \eqref{eq:ri.deltai.C}, the definition of $D_h$ in \eqref{def:damH}, and the fact that the function $\Xi(\cdot,\cdot)$ defined in \eqref{eq:Technical.varphi.def} is non-decreasing, and the last inequality is due to the definition of $\Xi(\cdot,\cdot)$ and the definition of $\Upsilon(\cdot )$  in \eqref{eq:tildezeta.Upsilon}.
\end{proof}

The proof of Lemma \ref{lem:mult.step} clearly uses the assumption that $\dom h$ is bounded (see (A1)). However, the above proof shows that the conclusion of Lemma \ref{lem:mult.step} also holds under the weaker assumption that the sequence $\{y^i\}$ generated by {\PAADMM} is bounded,
and $D_h$ and $\widehat \nabla_f$ are instead defined as $D_h:=\sup_i \|y^i-\bar{\mathrm{x}}\|$ and $\widehat \nabla_f:=\sup_{i} \| \nabla f(z^i)\|$, respectively.

\begin{lemma}[Multiplier and scaled residual bounds]
\label{lem:mult.bounds}
For every epoch $k$ generated within the $s$-th stage of {\PAADMM}, we have
\begin{align}
\|\tq^k\|
&\le
\max\{\|q^0\|,\Upsilon(C)\},
\label{eq:Lag.Mult.Bound}
\\
c\|Az^i-b\|
&\le
2\max\{\|q^0\|,\Upsilon(C)\},
\qquad
\forall i\in{\calI}_k(C),
\label{eq:Scaled.Residual.Bound}
\\
\|\tq^k-\tq^{k-1}\|
&\le
2\max\{\|q^0\|,\Upsilon(C)\},
\label{eq:Multiplier.Difference.Bound}
\end{align}
where $\Upsilon(C)$ and ${\calI}_k(C)$ are as in
\eqref{eq:tildezeta.Upsilon} and \eqref{def:I:k:C1}, respectively.
\end{lemma}

\begin{proof} 
Inequality \eqref{eq:Lag.Mult.Bound} follows by recursively using \eqref{step:induc:next:result} and the fact that $\tq^0=q^0$ due to  \eqref{def.variablesTilde}.  Inequality  \eqref{eq:Scaled.Residual.Bound} follows from \eqref{eq:FeasBound.proof.1}, the triangle inequality, \eqref{eq:Lag.Mult.Bound} with $k$ replaced by $k-1$, and the fact that $\tq^0=q^0$ due to  \eqref{def.variablesTilde}, i.e., 
\begin{align*}
c\|A z^i-b\|
\overset{\eqref{eq:FeasBound.proof.1}}\le \max\{\|\tilde q^{k-1}\|,\Upsilon(C)\}+\|\tq^{k-1}\|
\overset{\eqref{eq:Lag.Mult.Bound}}\leq 2\max\{\|\tilde q^{0}\|,\Upsilon(C)\}.
\end{align*}
Inequality \eqref{eq:Multiplier.Difference.Bound} follows from identity \eqref{eq:q.update}, the triangle inequality for norms, and \eqref{eq:Scaled.Residual.Bound} with $i=\cali_k$.
\end{proof}

The estimate \eqref{eq:Scaled.Residual.Bound} does not guarantee feasibility yet. Instead, it shows that the constraint violation at $z^i$ is about ${\cal O}(1/c)$. Later, the outer loop will turn this scaled residual estimate into the feasibility condition $\|A x^s-b\|\le \eta$ by increasing the penalty parameter.

\begin{lemma}\label{lem:mult.variation}
Let $k$ be an epoch index generated by {\PAADMM} during the $s$-th stage. Then
\begin{align}\label{eq.deltaq.bound}
\sum_{\ell=1}^{k-1}\|\tq^\ell-\tq^{\ell-1}\|^2\leq \frac{\tilde\zeta+c\zeta}{\alpha}\max_{1\leq\ell\leq k-1}\{\|\tq^{\ell}-\tq^{\ell-1}\|^2\},
\end{align}
where $\zeta$ and $\tilde \zeta$ are as in \eqref{eq:def.zeta2andblock_norm} and \eqref{eq:tildezeta.Upsilon}, respectively.
\end{lemma}
\begin{proof}
The case $k=1$ is trivial, since both sides of \eqref{eq.deltaq.bound} are zero. Hence, we assume that $k\geq 2$.
Using the fact that $\cali_{k-1}$ is the last index of ${\calI}_{k-1}$ and
noting the epoch termination criteria
in line \ref{ln:PAADMM:LMT} of {\PAADMM}, we see that
$\tilde T_{k-1}\leq \left[\sum_{j=1}^{\cali_{k-1}}(\|r^j\|^2+\delta_j)\right] /\left(\alpha R_{\cali_{k-1}}\right)$. The previous observation, \eqref{eq:C1b} with $i\in \{1,\ldots ,\cali_{k-1}\}$, the facts that $T_0=0$ by definition and
$T_{\cali_{k-1}}=\tilde T_{k-1}$ due to \eqref{def.variablesTilde}, imply that
\begin{align*}
\sum_{i=1}^{\cali_{k-1}}(\|r^i\|^2+\delta_i) \overset{\eqref{eq:C1b}}\leq \left(\tilde\zeta+c\zeta\right)\sum_{i=1}^{\cali_{k-1}}(T_{i}-T_{i-1}) =(\tilde\zeta+ c \zeta)\tilde T_{k-1}\leq \frac{(\tilde\zeta+ c \zeta)
\sum_{j=1}^{\cali_{k-1}}(\|r^j\|^2+\delta_j)}{\alpha R_{\cali_{k-1}}}.
\end{align*}

Rearranging the previous inequality, using the definition of $R_{i}$ as in line \ref{ln:PAADMM:R} of {\PAADMM} with $i=\cali_{k-1}$, and the fact that $\sum_{j=1}^{\cali_{k-1}}\|q^{j}-q^{j-1}\|^2=\sum_{\ell=1}^{k-1}\|\tq^{\ell}-\tq^{\ell-1}\|^2$ due to 
\eqref{eq:q-update} with $j=1,\ldots,\cali_{k-1}$, we conclude that \eqref{eq.deltaq.bound} holds.
\end{proof}

We now have all the ingredients needed to bound the potential sequence $\{T_i\}$ within a fixed stage. 
By combining Lemmas~\ref{lem:AL.gap}, ~\ref{lem:mult.bounds}, and~\ref{lem:mult.variation}, the next proposition derives a uniform upper bound on $T_i$ in terms of the problem parameters and the initial data of the current stage. It then uses the residual estimate in Lemma~\ref{lem:residual.potential} to obtain the corresponding bound on the number of inner iterations. 

\begin{prop}\label{prop:stage.complexity}
Let $s$ be a stage index of {\PAADMM} and define
\begin{align}\label{variable:inner.complexity}
 \Lambda(x^{s-1}, p^{s-1}, c_{s-1}\; ; C,\alpha) &:= 2(\phi_{\mathrm{sup}}-\phi_{\mathrm{inf}})
+ 2c_{s-1}\|Ax^{s-1}-b\|^2 \nonumber \\
&\qquad +\left[\frac{2}{c_{s-1}} + \frac{8}{\alpha}\left(\frac{\tilde \zeta}{c_{s-1}}+\zeta\right)\left(\frac{4B(1-\bar\tau)^2}{\bar \tau}+1\right)\right] \left(\|p^{s-1}\|^2+\Upsilon(C)^2\right),
\end{align}
where $(\alpha,C)$ is an input for {\PAADMM}, and $(\phi_{\mathrm{inf}},\phi_{\mathrm{sup}})$,  $\zeta$, $\tilde\zeta$, and $\Upsilon(C)$ are as in~\eqref{def:damH}, \eqref{eq:def.zeta2andblock_norm}, and \eqref{eq:tildezeta.Upsilon}, respectively.
Then, the following statements about the $s$-th stage of {\PAADMM} hold:
\begin{itemize}
    \item[(a)] it terminates after performing at most
    \begin{equation}\label{eq:inner.complexity}
1+ \left( \frac{\tilde\zeta+c_{s-1}\zeta}{\rho^2} \right) \Lambda(x^{s-1}, p^{s-1}, c_{s-1}\; ; C,\alpha)
\end{equation} 
inner iterations and the quintuple $(x^s,p^s, v^s, \varepsilon_s, c_s)$
    output by it  
   satisfies
\begin{equation}\label{eq:prop.inner.loop.inclusion}
v^s \in \nabla f(x^s) + \partial_{\varepsilon_s} h(x^s)+A^* p^{s}   \quad \text{and}\quad \|v^s\|^2+ \varepsilon_s \le \rho^2,
\end{equation}
and the inequalities
\begin{equation}\label{eq:prop.inner.loop.bounds}
c_{s}\|A x^{s}-b\| \le 4\max\{\|p^{s-1}\|,\Upsilon(C)\}  \quad \text{and} \quad
\quad \|p^s\| \le \max\{\|p^{s-1}\|,\Upsilon(C)\};
\end{equation}

\item[(b)] if $
   c_{s-1} \geq 2 [\max\{\|p^{s-1}\|,\Upsilon(C)\}]/\eta$,
then the quadruple $(x^s,p^s,v^s,\varepsilon_s)$
  is a $(\rho,\eta)$-stationary solution of problem  \eqref{initial.problem}-\eqref{eq:block_structure}. 
\end{itemize}
\end{prop}

\begin{proof}
(a) We divide the proof into two parts: we first prove \eqref{eq:inner.complexity}, and later \eqref{eq:prop.inner.loop.inclusion} and \eqref{eq:prop.inner.loop.bounds}. 

We first recall from line \ref{ln:PAADMM:stage.init} of {\PAADMM} that $(y^0,q^0,\lambda^0)=
    (x^{s-1},p^{s-1},\gamma^{s-1})$ and that $c_{s-1}=c$.
To prove \eqref{eq:inner.complexity}, we first claim that for every inner iteration $i$ generated within the $s$-th outer loop, we have that
\begin{equation}\label{eq.Ti.Lambda}
T_i \le \Lambda(y^{0}, q^{0}, c\; ; C,\alpha),
\end{equation}
where $\Lambda(y^{0}, q^{0}, c\; ; C,\alpha)=\Lambda(x^{s-1}, p^{s-1}, c_{s-1}\; ; C,\alpha)$ is as in \eqref{variable:inner.complexity}.
To prove the claim, let $k$ be the epoch index in which $i\in\calI_{k}$. 
It follows from  Lemma~\ref{lem:AL.gap} with $(x,p)=(\ty^{0};\tq^{0})=(y^0,q^0)$ (due to \eqref{def.variablesTilde}) and $(y,q)=(y^i;\tq^{k-1})$, that
\begin{align*}
{\cal L}_c(\ty^{0};\tq^{0}) - {\cal L}_c(y^i;\tq^{k-1})&\le \phi_{\mathrm{sup}}-\phi_{\mathrm{inf}}
+ c\|A\ty^{0}-b\|^2 +\frac{\|\tq^{0}\|^2+\|\tq^{k-1}\|^2}{2c}\\
&\le \phi_{\mathrm{sup}}-\phi_{\mathrm{inf}}
+ c\|Ay^0-b\|^2 +\frac{\|q^0\|^2+\Upsilon(C)^2}{c},
\end{align*}
where the second inequality is due to~\eqref{eq:Lag.Mult.Bound}.
Now, using \eqref{eq:Multiplier.Difference.Bound} and \eqref{eq.deltaq.bound}, we have
\[
\sum_{\ell=1}^{k-1}\|\tq^\ell-\tq^{\ell-1}\|^2\leq \frac{\tilde\zeta+c\zeta}{\alpha}\max_{1\leq\ell\leq k-1}\{\|\tq^{\ell}-\tq^{\ell-1}\|^2\}\leq       \frac{ 4(\tilde\zeta+c\zeta)}{\alpha}\left[\|q^0\|^2+\Upsilon(C)^2\right],
\]
and the claim now follows from \eqref{eq:Ti.bound.final}, the two inequalities above and \eqref{variable:inner.complexity}.

We now prove \eqref{eq:inner.complexity}.
Assume by contradiction that there exists an iteration index $i$ generated by {\PAADMM} such that
\begin{equation}\label{prel.contr:assum:i}
i > 1+ \left( \frac{\tilde\zeta+c\zeta}{\rho^2} \right) \Lambda(y^0, q^0, c\; ; C,\alpha).
\end{equation}
Since {\PAADMM} does not stop at any iteration index smaller than $i$,  the stopping criterion in line \ref{ln:PAADMM:stat.test} is violated at these iterations, i.e.,
$
\|r^j\|^2 +\delta_j > \rho^2
$ for every $j \le i-1$.
Hence, it follows from \eqref{eq:C1b}, the previous inequality, the fact that $T_0=0$ due to line \ref{ln:PAADMM:stage.init} of {\PAADMM}, and \eqref{eq.Ti.Lambda} that
\[
\frac{(i-1) \rho^2}{\tilde\zeta+c\zeta} < \frac{1}{\tilde\zeta+c\zeta}\sum_{j=1}^{i-1}(\|r^j\|^2+\delta_j) \overset{\eqref{eq:C1b}}\le \sum_{j=1}^{i-1}(T_j-T_{j-1}) = T_{i-1}-T_0 \le T_i \overset{\eqref{eq.Ti.Lambda}}{\leq}  \Lambda(y^0, q^0, c\; ; C,\alpha),
\]
which contradicts \eqref{prel.contr:assum:i}.
Thus, since $\Lambda(y^{0}, q^{0}, c\; ; C,\alpha)=\Lambda(x^{s-1}, p^{s-1}, c_{s-1}\; ; C,\alpha)$ we conclude that \eqref{eq:inner.complexity} holds.

We now prove the second part of (a), i.e., \eqref{eq:prop.inner.loop.inclusion} and \eqref{eq:prop.inner.loop.bounds}.
It follows from statement (a) that every stage generated by {\PAADMM} terminates. Hence, let $K$ be the last epoch generated during the $s$-th outer loop of {\PAADMM}. By line~\ref{ln:PAADMM:stage.output} of {\PAADMM}, the output of this stage satisfies
\[
(x^s,p^s,v^s,\varepsilon_s,\gamma^s)
=
(z^{\cali_K},q^{\cali_K},
r^{\cali_K},\delta_{\cali_K},\lambda^{\cali_K}).
\]
In particular, $p^s=q^{\cali_K}\in A(\R^n)$ in view of~\eqref{eq:q.in.Img.A.2}.
Using this identity, the inclusion in \eqref{eq:prop.inner.loop.inclusion} follows from Lemma~\ref{lem:residual.potential}, applied with $i={\cali_K}$. The residual bound in \eqref{eq:prop.inner.loop.inclusion} follows from the stopping condition in line~\ref{ln:PAADMM:stat.test} of {\PAADMM}. Finally, the bounds in \eqref{eq:prop.inner.loop.bounds} follow from Lemma~\ref{lem:mult.bounds} with $i=\cali_K$ and the facts that $q^0=p^{s-1}$ and $\tilde q^K=q^{\cali_K}$.

(b) It follows from the first inequality in \eqref{eq:prop.inner.loop.bounds} that, if $c_{s-1} \ge [2\max\{\|p^{s-1}\|,\Upsilon(C)\}]/\eta$,
then $\|Ax^s-b\|\le \eta$. Combining this feasibility bound with the \textit{approximate stationarity inclusion} and \textit{stationarity residual bound} in \eqref{eq:prop.inner.loop.inclusion}, we conclude that the quadruple $(x^s,p^s,v^s,\varepsilon_s)$
is a $(\rho,\eta)$-stationary solution of \eqref{initial.problem}--\eqref{eq:block_structure}.
\end{proof}

\vspace{1em}

It follows from Proposition~\ref{prop:stage.complexity}(a) that the $s$-th stage returns in its line~\ref{ln:PAADMM:stage.output} a quadruple $(x^s, p^s, v^s, \varepsilon_s)$ which is a $\rho$-stationary solution (see \eqref{eq:rho.eta.solution} and the paragraph immediately following it), but $x$ is not necessarily $\eta$-feasible. However, Proposition~\ref{prop:stage.complexity}(b) guarantees that, if $c_{s-1}$ is already large enough, i.e., $c_{s-1} = \Omega(\eta^{-1})$,
then $x$ is $\eta$-feasible,
and hence that $(x^s, p^s, v^s, \varepsilon_s)$ is a $(\rho,\eta)$-stationary solution of \eqref{initial.problem}-\eqref{eq:block_structure}.

\vspace{1em}

\noindent
{\bf Proof of Theorem~\ref{thm:main}.}
(a) Both inequalities in \eqref{eq:prop:input.bounds} follow from \eqref{eq:prop.inner.loop.bounds}, a simple induction argument, and the fact that $p^0=0$ due to line \ref{ln:PAADMM:init} of {\PAADMM}.

(b) It follows from Proposition~\ref{prop:stage.complexity} that the total number of inner iterations performed by {\PAADMM} during the $s$-th outer loop is bounded.
Hence, it suffices to show that 
\[
\Lambda(x^{s-1},p^{s-1},c_{s-1};C,\alpha)\leq \Gamma(x^0;C,\alpha),
\] 
where $\Gamma(x^0;C,\alpha)$ is as in \eqref{eq:def.Gamma.complexity}.

Before showing the above inequality, we first show that, for every stage index $s$,
\[
c_{s-1}\|Ax^{s-1}-b\|^2\leq \frac{16\Upsilon^2(C)}{c_0}.  %\quad \text{for every stage index} \quad s.
\] 
Indeed, for $s=1$, the definition of $c_0$ in line~\ref{ln:PAADMM:init} of {\PAADMM} implies that $\|Ax^0-b\|\le 1+\|Ax^0-b\|=1/c_0$, and hence
\[
c_0\|Ax^0-b\|^2\le \frac{1}{c_0}\le \frac{16\Upsilon^2(C)}{c_0},
\]
where the last inequality follows from the fact that $\Upsilon(C)\geq 1$ due to its definition in~\eqref{eq:tildezeta.Upsilon}. Moreover, the first inequality in~\eqref{eq:prop.inner.loop.bounds} and the fact that $c_{s-1}\geq c_0$ show that the inequality also holds for
$s>1$, and thus it holds for any $s \ge 1$.

Using the last conclusion, the definition of $\Lambda(x^{s-1},p^{s-1},c_{s-1};C,\alpha)$ as in \eqref{variable:inner.complexity}, 
the fact that $c_{s-1}\ge c_0$, and that $\|p^{s-1}\|\le \Upsilon(C)$ due to (a), we have
\begin{align*}
& \Lambda(x^{s-1},p^{s-1},c_{s-1};C,\alpha) \\
& \leq 2(\phi_{\mathrm{sup}}-\phi_{\mathrm{inf}})
+ 2c_{s-1}\|Ax^{s-1}-b\|^2 +\left[\frac{2}{c_{0}} + \frac{8}{\alpha}\left(\frac{\tilde \zeta}{c_{0}}+\zeta\right)\left(\frac{4B(1-\bar\tau)^2}{\bar\tau}+1\right)\right] \left(\|p^{s-1}\|^2+\Upsilon(C)^2\right)  \\
&\leq  2(\phi_{\mathrm{sup}}-\phi_{\mathrm{inf}}) + \frac{32\Upsilon(C)^2}{c_{0}} + \left[\frac{4}{c_{0}} + \frac{16}{\alpha}\left(\frac{\tilde \zeta}{c_{0}}+\zeta\right)\left(\frac{4B(1-\bar\tau)^2}{\bar\tau}+1\right)\right] \Upsilon(C)^2 \\
&=2(\phi_{\mathrm{sup}}-\phi_{\mathrm{inf}}) +(1+\|Ax^0-b\|)\Upsilon(C)^2\left[36+ \frac{16\tilde \zeta }{\alpha}\left(\frac{4B(1-\bar\tau)^2}{\bar\tau}+1\right)\right] + \frac{16\zeta \Upsilon(C)^2}{\alpha}\left(\frac{4B(1-\bar\tau)^2}{\bar\tau}+1\right)\\
&=\Gamma(x^0;C,\alpha),
\end{align*}
which concludes the proof.

(c) Assume, for the sake of contradiction, that {\PAADMM} generates a stage $s$ such that $s>1+\log\left(1+\frac{4\Upsilon(C)}{c_0\eta}\right)$, and hence satisfying
\[
c_{s-1}=2^{s-1}c_0
>
\frac{4\Upsilon(C)}{\eta}.
\]
Note that $s\ge2$ and $\|p^{s-2}\|\le\Upsilon(C)$, where the latter follows from Theorem~\ref{thm:main}(a) applied to the $(s-2)$-th stage if $s\ge3$, and from $p^0=0$ if $s=2$. Hence,
\[
c_{s-2}=\frac{c_{s-1}}{2}>\frac{2\Upsilon(C)}{\eta}=\frac{2\max\{\|p^{s-2}\|,\Upsilon(C)\}}{\eta},
\]
and Proposition~\ref{prop:stage.complexity}(b) applied to the $(s-1)$-th stage implies that the quadruple $(x^{s-1},p^{s-1},v^{s-1},\varepsilon_{s-1})$ is a $(\rho,\eta)$-stationary solution of \eqref{initial.problem}--\eqref{eq:block_structure}. Therefore, {\PAADMM} would have terminated at the end of the $(s-1)$-th stage, contradicting the assumption that stage $s$ was generated \qedflush 

\section{Numerical Experiments}\label{sec:numerics}

This section presents the numerical performance of {\PAADMM} on a nonconvex linearly constrained quadratic problem. We first explain in Subsection~\ref{sub:numerics.setup} the experimental setup, and in Subsection~\ref{sub:numerics.QP} we show the results of three methods in the quadratic problem.

\subsection{Experimental Setup}\label{sub:numerics.setup}

The computational times presented for the parallel variants are
idealized parallel execution times. Specifically, during each
{\PABIPP} call, the block subproblems are solved and timed individually, and
the time attributed to the parallel block-solution phase corresponds to the maximum of
the individual block times. All computations that depend on the outputs of
each block, such as the construction of the relaxed iterate and the
subsequent residual calculations, are included in full.

The convergence analysis of {\PAADMM} uses the relaxation parameter $\bar\tau=1/[4(B-1)]$ (see line~\ref{ln:PABIPP:tau} of {\PABIPP}).
Although this choice guarantees the required descent property (see Lemma~\ref{lem:suff.decrease}), it can be
conservative in practice, especially as the number of blocks increases. We therefore use an adaptive search method for improving the value of $\bar \tau$. The complete search procedure is described
in Appendix~\ref{app:adaptive.tau}.

We compare {\PAADMM} with two other methods. The first is the cyclic adaptive method {\AADMM} of~\cite{MaiaEtAl2026}, i.e., the Gauss--Seidel counterpart of {\PAADMM}. The second, denoted by {\JAADMM}, is {\PAADMM} with its Lagrange multiplier test replaced by the one of {\AADMM}~\cite{MaiaEtAl2026}. We also attempted to run two other variants of {\PAADMM}. The first is a penalty variant that never updates the Lagrange multiplier, i.e., the method without lines 12-19 of {\PAADMM} and the Lagrange multiplier updates in lines 6 and 9. The second is a vanilla ADMM variant that skips the LMT in line~\ref{ln:PAADMM:LMT} and updates the Lagrange multiplier at every inner iteration, as in the standard ADMM. Neither of these two variants converged on the tested instances. Finally, we attempted to run the dampened proximal ADMM (DP-ADMM) of~\cite{KongMonteiro2024} with many choices of its parameter pair $(\chi,\theta)$, but none of the tested choices led it to satisfy the stopping criterion below within the iteration limit. We therefore omit these three methods from the tables.

For each method, the tables report the total number of iterations
(``Iter''), the computational time (``Time''), and the objective value at
the returned point (``$\phi$''). We also report ``Mults'', the total number of (accepted) Lagrange multiplier updates.

The blocks of the generated instances of
\eqref{initial.problem}--\eqref{eq:block_structure} are taken to have the
same dimension,
\[
    n_1=\cdots=n_B=\mathsf{n},
    \qquad
    n=B\mathsf{n}.
\]
Hence, apart from problem-specific parameters, an instance is characterized
by the triple $(B,\mathsf{n},\mathsf{m})$,
where $B$ is the number of blocks, $\mathsf{n}$ is the block dimension, and
$\mathsf{m}$ is the number of rows of the constraint matrix $A$. The values of $(\omega,B,\mathsf{n},\mathsf{m})$ used in each experiment are reported
at the beginning of the corresponding table. Each block matrix
$A_t\in\mathbb{R}^{\mathsf{m}\times\mathsf{n}}$ is generated with i.i.d.\
standard-normal entries. The initial penalty parameter and multiplier are
chosen according to line~\ref{ln:PAADMM:init} of {\PAADMM}, namely, $c_0=1/\bigl(1+\|Ax^0-b\|\bigr)$ and $p^0=0$. We further set $C=10^3\rho\bigl(1+\|\nabla f(x^0)\|\bigr)$,
and initialize every adaptive block stepsize with the common value $\lambda_t^0=10$ for $t=1,\ldots,B$.

In the two experiments, the block subproblem
\eqref{algo.PABIPP-subproblem} is solved using the ADAP-FISTA routine from \cite{sujanani2023adaptive} with $(M_0,\beta,\mu,\chi)=(1,1.2,0.5,0.001)$.
The routine is terminated once it produces a $(1/32,y_t)$-stationary point.
Whenever the required conditions for accepting the current block solution
are not satisfied, the current stepsize $\lambda_t$ is halved and {ADAP-FISTA} is restarted.

Each method is allowed a maximum of $500{,}000$ iterations. Otherwise, it
terminates upon producing an approximate stationary quadruple
$(x^{\splus},p^{\splus},v^{\splus}, \varepsilon^{\splus})$ satisfying
\[
    v^{\splus}
    \in
    \nabla f(x^{\splus})+\partial_{\varepsilon^{\splus}} h(x^{\splus})+A^*p^{\splus},
\]
together with the relative stationarity and feasibility conditions
\[
    \frac{\|v^{\splus}\|}
         {1+\|\nabla f(x^0)\|}
    \leq \rho,
    \qquad
    \frac{\|Ax^{\splus}-b\|}
         {1+\|Ax^0-b\|}
    \leq \eta,
    \qquad
    \rho=\eta=10^{-5}.
\]

All experiments were implemented in MATLAB R2024b and executed on a
macOS machine equipped with an Apple M5 processor (10 cores) and
32~GB of memory. For brevity, the experiments below focus on randomly generated dense instances.

\subsection{Nonconvex Distributed Quadratic Programming}
\label{sub:numerics.QP}

We consider a class of nonconvex distributed quadratic programming
problems (DQP) with $B$ blocks. This problem is useful for testing the
parallel variants since both the objective function and the linear
constraint are naturally decomposed according to the block structure.

Given $(\mathsf{m},\mathsf{n})\in\mathbb{N}_{++}^2$ satisfying
$\mathsf{m}<B\mathsf{n}$, we consider
\begin{align}
\label{eq:num.problem1}
\min_{x=(x_1,\ldots,x_B)\in\mathbb{R}^{B\mathsf{n}}}
\quad &
\sum_{i=1}^B
\left[
    \frac12 \langle x_i,P_i x_i\rangle
    +\langle x_i,r_i\rangle
\right]
\nonumber\\
\text{s.t.}\quad &
\|x_i\|_1\leq \omega,
\qquad i=1,\ldots,B,\\
&
\sum_{i=1}^B A_i x_i=b,
 \nonumber \end{align}
where $\omega>0$,
$P_i\in\mathbb{R}^{\mathsf{n}\times\mathsf{n}}$ is symmetric and
indefinite,
$r_i\in\mathbb{R}^{\mathsf{n}}$,
$A_i\in\mathbb{R}^{\mathsf{m}\times\mathsf{n}}$,
and $b\in\mathbb{R}^{\mathsf{m}}$. Problem~\eqref{eq:num.problem1} is a special case of
\eqref{initial.problem}--\eqref{eq:block_structure} with
\[
f(x)
=
\sum_{i=1}^B
\left[
    \frac12\langle x_i,P_i x_i\rangle
    +\langle x_i,r_i\rangle
\right]
\quad \text{and} \quad h_i(x_i)
=
\delta_{\{u\in\mathbb{R}^{\mathsf n}:\|u\|_1\leq\omega\}}(x_i),
\qquad i=1,\ldots,B.
\]

We generate the instances as follows. First, after generating the matrices
$A_i$, we construct a feasible vector
$x^b=(x_1^b,\ldots,x_B^b)$ by drawing each block randomly and scaling it, if
necessary, so that $\|x_i^b\|_1 < \omega$. Next, we set $b=\sum_{i=1}^B A_i x_i^b$. This ensures that the generated instance is feasible. We generate the starting point $x^0$ independently in the same way, with $\|x_i^0\|_1<\omega$ for every block. To build the quadratic terms, for each block $i$ we first generate an orthonormal matrix $Q_i$. Then we create a diagonal matrix $D_i$ where one third of the diagonal entries are set to zero, and the rest are sampled uniformly from $[-10,10]$. We also require that $D_i$ has at least one negative entry. Finally, we set $P_i=Q_i^\top D_iQ_i$. This gives a symmetric indefinite matrix $P_i$, which leads to a nonconvex quadratic objective. More specifically, if $m_i<0$ denotes the smallest diagonal entry of $D_i$, then $\lambda_{\min}(P_i)=m_i$, and each block function
$x_i
\mapsto
\frac12\langle x_i,P_i x_i\rangle+\langle x_i,r_i\rangle
$
is $|m_i|$-weakly convex. Finally, the entries of each vector $r_i$ are
drawn independently from the standard normal distribution.

The numerical results for this problem class are reported in
Table~\ref{tab:DQP}.

\begin{table}[ht]
\centering
\scriptsize
\begin{adjustbox}{max width=\textwidth}
\begin{tabular}{c c c | r r r | r r r | r r r}
\hline
\multicolumn{3}{c|}{Instance} & \multicolumn{3}{c|}{\PAADMM} & \multicolumn{3}{c|}{\AADMM} & \multicolumn{3}{c}{\JAADMM} \\
$\omega$ & $B$ & $(m,n)$ & Time & Iters/Mult & $\phi$ & Time & Iters/Mult & $\phi$ & Time & Iters/Mult & $\phi$ \\
\cmidrule(lr){1-3}\cmidrule(lr){4-6}\cmidrule(lr){7-9}\cmidrule(lr){10-12}
100 & 16 & (32,32) & \textbf{25.978} & \textbf{59916}/2067 & -1.099e+05 & 72.599 & 128550/11 & -1.004e+05 & 40.102 & 128017/10 & -1.136e+05 \\
 &  & (64,64) & 84.089 & 133061/421 & -4.804e+04 & \textbf{46.123} & \textbf{16019}/9 & -4.534e+04 & 59.850 & 99916/10 & -4.919e+04 \\
 &  & (64,96) & 207.265 & 306514/127 & -5.068e+04 & 111.872 & \textbf{47308}/10 & -4.871e+04 & \textbf{32.853} & 48127/9 & -5.210e+04 \\
 &  & (96,64) & \textbf{77.157} & 115383/825 & -5.638e+04 & 152.699 & \textbf{91877}/11 & -5.313e+04 & 158.496 & 250393/11 & -5.610e+04 \\
\\[0.35em]
 & 32 & (32,32) & \textbf{28.865} & \textbf{39375}/32 & -2.001e+05 & 175.949 & 164027/9 & -1.898e+05 & 223.829 & 317578/9 & -2.025e+05 \\
 &  & (64,64) & \textbf{35.101} & 41194/64 & -1.362e+05 & 151.852 & \textbf{40996}/10 & -1.382e+05 & 264.708 & 322900/10 & -1.433e+05 \\
 &  & (64,96) & 399.247 & 369175/47 & -8.348e+04 & \textbf{369.075} & \textbf{137465}/9 & -8.771e+04 & 392.666 & 380323/10 & -8.676e+04 \\
 &  & (96,64) & 224.996 & 236909/103 & -1.151e+05 & * & * & * & \textbf{126.373} & \textbf{144487}/10 & -1.124e+05 \\
\\[0.35em]
 & 64 & (32,32) & \textbf{27.700} & \textbf{23603}/7 & -3.313e+05 & 122.846 & 50943/7 & -3.082e+05 & 90.178 & 83737/7 & -3.343e+05 \\
 &  & (64,64) & \textbf{108.986} & \textbf{73533}/8 & -2.502e+05 & 357.803 & 94601/8 & -2.444e+05 & 375.691 & 260469/8 & -2.523e+05 \\
 &  & (64,96) & \textbf{33.314} & \textbf{19646}/8 & -2.280e+05 & * & * & * & 189.575 & 114267/8 & -2.311e+05 \\
 &  & (96,64) & \textbf{83.250} & \textbf{52285}/10 & -2.796e+05 & * & * & * & 440.221 & 286355/9 & -2.795e+05 \\
\\[0.35em]
 & 128 & (32,32) & \textbf{21.994} & \textbf{10261}/6 & -7.584e+05 & * & * & * & 159.819 & 78146/7 & -7.570e+05 \\
 &  & (64,64) & \textbf{22.334} & \textbf{11074}/7 & -5.673e+05 & 1301.873 & 211845/7 & -5.116e+05 & 888.601 & 473254/7 & -5.647e+05 \\
 &  & (64,96) & \textbf{146.600} & \textbf{42737}/8 & -4.889e+05 & 1101.274 & 106021/7 & -4.523e+05 & 509.547 & 152617/8 & -4.953e+05 \\
 &  & (96,64) & \textbf{452.062} & \textbf{155410}/9 & -5.748e+05 & * & * & * & * & * & * \\
\hline
\end{tabular}
\end{adjustbox}
\caption{Performance of the ADMM variants for the DQP problem.}
\label{tab:DQP}
\end{table}

The results in Table~\ref{tab:DQP} show that {\PAADMM} benefits substantially from the parallel block updates. First, {\PAADMM} converged in all $16$ instances, whereas {\AADMM} and {\JAADMM} did not converge in $5$ and $1$ instances, respectively (marked with * in Table~\ref{tab:DQP}). Second, {\PAADMM} achieved the shortest time in $12$ of the $16$ instances. The largest improvement occurs for $B=128$ and $(\mathsf{m},\mathsf{n})=(64,64)$, where {\PAADMM} is approximately $58$ times faster than {\AADMM} and approximately $40$ times faster than {\JAADMM}.

Third, {\PAADMM} also performs well in terms of iterations, obtaining the fewest iterations in $10$ of the $16$ instances. Since {\PAADMM} and {\JAADMM} differ only in the Lagrange multiplier test, their comparison highlights the benefit of the new LMT. Compared to {\JAADMM}, {\PAADMM} required fewer iterations in $12$ of the $15$ instances in which both converged.

\section{Concluding Remarks}\label{sec:conclusion}

This paper suggests several directions for further investigation, particularly extending the proposed framework to more general objective functions and block-update schemes. An important extension is to consider a nonseparable differentiable component $f$, as in \cite{MaiaEtAl2026}. Another direction is to investigate block-update schemes that interpolate between the fully parallel method considered here and the cyclic method of \cite{MaiaEtAl2026}. One such scheme partitions the blocks into groups, updates the blocks within each group in parallel, and processes the groups sequentially. It would also be worthwhile to weaken some assumptions on the nonsmooth component $h$. In particular, our analysis assumes that $\dom h$ is bounded and that $h$ is Lipschitz continuous on its domain (see (A1)). 

Finally, the proposed framework also yields an adaptive parallel proximal penalty method. Indeed, choosing the input parameter $C$ of {\PAADMM} such that $C>\rho$ ensures that the test in line~\ref{ln:PAADMM:LMT} is never satisfied, so the LMT occurs only at the end of each stage. Since these updates are not needed for the convergence analysis, we can omit them altogether. The resulting method is an adaptive parallel proximal penalty scheme satisfying convergence guarantees.

\appendix

\section{Technical Details}\label{app:technical}

This section presents technical lemmas used in the iteration-complexity analysis. Lemma~\ref{lem:AL.block.identity} is used to write the decrease of $\Lc$ from the parallel update as a telescoping sum of single-block decreases. Lemmas~\ref{lem:eps.subgrad} and~\ref{lem:mult.bound.aux} are taken from~\cite{MaiaEtAl2026}. Lemma~\ref{lem:eps.subgrad} provides a lower bound for a strongly convex function using an inexact subgradient. Lemma~\ref{lem:mult.bound.aux} bounds the current Lagrange multiplier in terms of the previous multiplier and the residual pair.

\begin{lemma}
\label{lem:AL.block.identity}
For given $t\in\{1,\ldots,B\}$, $x,y\in \dom h$, and $(p,c)\in A(\R^n)\times\R_{++}$, it holds that
\begin{align}\label{eq:lemma.renato.paper.equation}
    \Lc(y_{<t},x_t,x_{> t};p) - \Lc(y_{<t}, y_{t},x_{>t};p) & = \Lc(x;p) - \Lc(x_{<t},y_{t},x_{>t};p) \nonumber\\
    &\qquad - c\sum_{i=1}^{t-1}\inner{A_t(y_t-x_t)}{A_i(y_i-x_i)}
\end{align}
\end{lemma}
\begin{proof} 
The case $t=1$ is immediate. We assume now that $t\ge 2$, and  to simplify notation we define the function $\mathrm{L}:{\cal H}_1\times\cdots\times{\cal H}_{t-1}\to\mathbb{R}$ by
\begin{equation*}
    \vecu_{<t} \mapsto \mathrm{L}(\vecu_{<t}):=\Lc(\vecu_{<t},x_{t},x_{>t};p) - \Lc(\vecu_{<t},y_{t},x_{>t};p).
\end{equation*}
With this definition, \eqref{eq:lemma.renato.paper.equation} is equivalent to showing that
\begin{equation}\label{eq:proof.renato.paper.goal}
    \mathrm{L}(y_{<t})
    = \mathrm{L}(x_{<t})
    - c\sum_{i=1}^{t-1}\inner{A_t(y_t-x_t)}{A_i(y_i-x_i)}.
\end{equation}
We now derive \eqref{eq:proof.renato.paper.goal}. Using \eqref{eq:block_structure} and \eqref{DP:AL:F}, we see that the $\vecu_{<t}$-dependent parts of the terms of $\mathrm{L}$ involving $f$, $h$, and $p$ cancel out. Hence, $\mathrm{L}$ is the restriction to ${\cal H}_1\times\cdots\times{\cal H}_{t-1}$ of the quadratic function
\begin{equation*}
    \vecu_{<t}\mapsto \kappa+\frac{c}{2}\|A(\vecu_{<t},x_{t},x_{>t})-b\|^2-\frac{c}{2}\|A(\vecu_{<t},y_{t},x_{>t})-b\|^2,
\end{equation*}
where $\kappa\in\R$ does not depend on $\vecu_{<t}$. The gradient of this quadratic function is
\begin{align*}
     &c[A_1\ldots A_{t-1}]^*\big(A(\vecu_{<t},x_{t},x_{>t})-b\big)-c[A_1\ldots A_{t-1}]^*\big(A(\vecu_{<t},y_{t},x_{>t})-b\big)\\
     &\qquad=c[A_1\ldots A_{t-1}]^*A_t(x_{t}-y_{t}),
\end{align*}
which does not depend on $\vecu_{<t}$, and therefore this quadratic function is affine. Consequently,
\begin{align*}
\text{L}({y}_{<t}) - \text{L}(x_{<t}) &= \inner{c[A_1\ldots A_{t-1}]^*A_t(x_{t}-y_{t})}{y_{<t} - x_{<t}} \\
    &= c\inner{A_t(x_{t}-y_{t})}{[A_1\ldots A_{t-1}](y_{<t} - x_{<t})} \\
    &=-c\sum_{i=1}^{t-1}\inner{A_t(y_t-x_{t})}{A_i(y_i-x_i)},
\end{align*}
which shows that~\eqref{eq:proof.renato.paper.goal} holds true, and we conclude the proof.
\end{proof}

\begin{lemma}\label{lem:eps.subgrad}
		Assume that $\xi>0$, $\psi \in \bConv{n}$ and positive definite real-valued $n \times n$-matrix $Q$ are such that
		$\psi - (\xi/2) \|\cdot\|^2_Q$ is convex and let
		$(y,v,\eta) \in \R^n \times \R^n \times  \R_{\splus}$ be such that $v\in \partial_\eta \psi(y)$.
		Then, %for any $ u \in \dom \psi$,
		for any $\varsigma>0$,
			\begin{equation*}
		\psi(u) \ge \psi(y) + \Inner{v}{u-y} - (1+\varsigma^{-1})\eta+ \frac{(1+\varsigma)^{-1}\xi}{2} \|u-y\|_Q^2  \quad \forall u \in \R^n.
		\end{equation*}
	\end{lemma}

\begin{proof}
    See \cite[Lemma A.4]{MaiaEtAl2026}.
\end{proof}

\begin{lemma}\label{lem:mult.bound.aux}
Assume that
$b \in \R^{l}$,
linear operator $A:\mathbb{R}^n \to \mathbb{R}^l$
and function $h(\cdot)$,
satisfy assumptions (A3)
and (A1), respectively.  
If $(q^-,\varrho) \in A(\R^n) \times (0,\infty)$ and $(z,q,r,\delta) \in  \dom h \times \mathbb{R}^l \times \mathbb{R}^n\times \R_+$ satisfy
\begin{equation}\label{eq:cond:Lem:A1}
q=q^-+\varrho(Az-b)\quad\text{and}\quad  r \in \partial_\delta h(z)+A^{*}q 
, %\textrm{[choose another letter for $\chi$]}
\end{equation}
then we have 
\begin{equation}\label{q-bound-2}
 \|q\|\leq \max\left\{\|q^-\|, \ \Xi\Big(\|z-\bar x\|\ , \ \|r\|+\delta \Big) \right\}
\end{equation}
where  $\mathrm{\bar x}$ is as in (A3),
\begin{equation}\label{eq:Technical.varphi.def}
\Xi(s,t) := 
%\upsilon(t) := 
\frac{t+ \left(s+\bar d \, \right)(M_h+t)}{\bar d \nu^{+}_A} \quad \forall (s,t)
\in \R_+^2,
\end{equation}
$M_h$ and $\bar d>0$ are as in (A1) and (A3), respectively, and $\nu^+_A$ is  the smallest positive singular value of $A$.
\end{lemma}

\begin{proof}
    See \cite[Lemma A.3]{MaiaEtAl2026}.
\end{proof}

\section{Adaptive Relaxation Parameter}\label{app:adaptive.tau}

This section describes the variant of {\PABIPP} used in the numerical experiments of Section~\ref{sec:numerics}, in which the relaxation parameter is chosen adaptively. 

{\floatname{algorithm}{Subroutine}
\renewcommand{\thealgorithm}{}
\begin{algorithm}[H]
\setstretch{1}
\caption{{\PABIPP} with adaptive relaxation parameter}
\begin{algorithmic}[1]
\Statex \hskip-1.8em
\Require $(y,q,\lam,c)
\in {\cal H}\times A(\R^n)\times\R^B_{++}\times\R_{++}$
\Ensure $(y^{\splus},\lam^{\splus};z^{\splus},r^{\splus},\delta^{\splus})
\in {\cal H}\times\R^B_{++}\times{\cal H}\times\R^n\times\R_{+}$

\vspace{1em}

\Statex \hskip-2.8em \textbf{\textcolor{algemph}{STEP 1: Parallel adaptive block IPP iteration}}

\Statex \hskip-2.8em
Same as STEP 1 of {\PABIPP}. Output: $(z^{\splus}, \lambda^{\splus})$

\Statex

\Statex \hskip-2.8em
\textbf{\textcolor{algemph}{STEP 2: Bisection search for the relaxation parameter $\bar\tau^{\splus}$}}

\State $\tau_1 \gets \bar\tau$ and $\tau_2 \gets 1$
\For{$i=1,\ldots,8$}
    \State $\tau^{\mathrm{mid}} \gets (\tau_1+\tau_2)/2$
    \If{
    \[
    \hspace{1.5cm} \Lc(y;q) - \Lc(y(\tau^{\mathrm{mid}});q)
    \ge
    \frac{(\tau^{\mathrm{mid}})^2}{2}
    \sum_{t=1}^B \frac{\|z_t^{\splus}-y_t\|^2}{\lambda_t^{\splus}}
    +
    \frac{c\tau^{\mathrm{mid}}}{8}
    \sum_{t=1}^B \|A_t(z_t^{\splus}-y_t)\|^2
    \]
    }
        \State $\tau_1 \gets \tau^{\mathrm{mid}}$ \Comment{midpoint is feasible}
    \Else
        \State $\tau_2 \gets \tau^{\mathrm{mid}}$ \Comment{midpoint is infeasible}
    \EndIf
\EndFor
\State $\bar \tau^{\splus} \gets \tau_1$ and $y^{\splus}\gets y+\bar\tau^{\splus}(z^{\splus}-y)$

\vspace{1em}

\Statex \textbf{\textcolor{algemph}{STEP 3: Computation of the residual pair $(r^{\splus},\delta^{\splus})$}}

%\Statex \hskip-2.8em
\Statex Same as STEP 2 of {\PABIPP}

\end{algorithmic}
\end{algorithm}
\addtocounter{algorithm}{-1}
\renewcommand{\thealgorithm}{\arabic{algorithm}}
}

The value $\bar\tau=1/[4(B-1)]$ in line~\ref{ln:PABIPP:tau} of {\PABIPP} guarantees the sufficient decrease \eqref{eq:renato.paper3}, but it can be conservative when $B$ is large. The subroutine above therefore inserts, between STEPs~1 and~2 of {\PABIPP}, a bisection search for a larger relaxation parameter. This search is carried out after the block subproblems have been solved, so $(z^{\splus},\lambda^{\splus})$ is as in line~\ref{ln:PABIPP:update} of {\PABIPP} and $y(\tau)=y+\tau(z^{\splus}-y)$ is as in \eqref{def:z_tau}. A trial value $\tau$ is declared admissible if it satisfies the test in STEP~2, i.e., \eqref{eq:renato.paper3} with $\bar\tau$ replaced by $\tau$. Starting from the interval $[\bar\tau,1]$, the search performs eight bisection steps. Since $\bar\tau$ is admissible (see Lemma~\ref{lem:suff.decrease}), $\tau_1$ is admissible throughout the search and is the largest admissible value tested so far. Hence, the output $\bar\tau^{\splus}\in[\bar\tau,1]$ is admissible, and it is used to construct $y^{\splus}=y(\bar\tau^{\splus})$. Moreover, since the right-hand side of \eqref{eq:renato.paper3} is increasing in $\tau$, $y^{\splus}$ satisfies \eqref{eq:adapt.Delta.Lc.ineq}, which is the only property of $y^{\splus}$ used in the proof of Proposition~\ref{prop:PABIPP}. Hence, the conclusions of Proposition~\ref{prop:PABIPP} remain valid for this variant.

\bibliographystyle{plainurl}   % or abbrv, ieeetr, alpha, etc.
\bibliography{Jacobi_ADMM_references} 

\end{document}